\documentclass[12pt, leqno]{amsart}

\usepackage{graphicx}
\usepackage{amssymb,amsmath,amsthm,amscd,accents}
\usepackage{mathrsfs}
\usepackage{lscape}
\usepackage{enumerate}
\usepackage{upgreek}
\usepackage{xspace}
\usepackage[usenames,dvipsnames]{color}
\usepackage[colorlinks=true, pdfstartview=FitV,
 linkcolor=blue,citecolor=blue,urlcolor=blue]{hyperref}
\usepackage{tikz}
\usetikzlibrary{arrows.meta,arrows}
\usepackage{bm}
\usepackage{marginnote}
\usepackage{verbatim}
\usepackage{comment}
\usepackage[normalem]{ulem}  

\usepackage[all]{xy}

\allowdisplaybreaks[3]

\newcommand{\arxiv}[1]{\href{http://arxiv.org/abs/#1}{\texttt{arXiv:#1}}}

\newcommand{\nc}{\newcommand}
\numberwithin{equation}{section}

\newcommand{\ql}{quasi-Laurent\xspace}

\newcommand{\dsc}{independent\xspace}

\nc{\berm}[1]{\begin{red}{}\marginnote{\fbox{\scshape\lowercase{M}}}%
#1}

\newcommand{\bero}[1]{\begin{red}{}\marginnote{\fbox{\scshape\lowercase{O}}}%
#1}

\newcommand{\berE}[1]{\begin{red}{}\marginnote{\fbox{\scshape\lowercase{E}}}%
#1}

\newcommand{\berMH}[1]{\begin{red}{}\marginnote{\fbox{\scshape\lowercase{MH}}}%
#1}

\nc{\hs}{\hspace*}
\nc{\ms}{\mspace}

\nc{\qR}[1]{q_{\mspace{-1mu}\raisebox{-.5ex}{${\scriptstyle{#1}}$}}}
\nc{\bfgi}{\bfg_{\mspace{-1mu}\raisebox{-.5ex}{${\scriptstyle{\ii}}$}}}

\theoremstyle{plain}
\newtheorem{lemma}{Lemma}[section]
\newtheorem{proposition}[lemma]{Proposition}
\newtheorem{theorem}[lemma]{Theorem}

\newtheorem{corollary}[lemma]{Corollary}

\theoremstyle{definition}
\newtheorem{remark}[lemma]{Remark}
\newtheorem{example}[lemma]{Example}

\newtheorem{definition}[lemma]{Definition}

\newtheorem*{convt}{Convention}

\nc{\Cor}{\begin{corollary}}
\nc{\encor}{\end{corollary}}
\nc{\Rem}{\begin{remark}}
\nc{\enrem}{\end{remark}}
\nc{\Lemma}{\begin{lemma}}
\nc{\enlemma}{\end{lemma}}
\nc{\Prop}{\begin{proposition}}
\nc{\enprop}{\end{proposition}}

\newcommand{\Lto}{\longrightarrow}

\renewcommand{\le}{\leqslant}
\renewcommand{\ge}{\geqslant}
\renewcommand{\preceq}{\preccurlyeq}

\newcommand{\head}{{\operatorname{hd}}}

\newcommand{\hconv}{\mathbin{\scalebox{.9}{$\nabla$}}}
\newcommand{\sconv}{\mathbin{\scalebox{.9}{$\Delta$}}}

\newcommand{\seteq}{\mathbin{:=}}
\newcommand{\conv}{\mathop{\mathbin{\mbox{\large $\circ$}}}}
\newcommand{\soplus}{\mathop{\mbox{\normalsize$\bigoplus$}}\limits}

\newcommand{\tens}{\mathop\otimes}

\newcommand{\gmod}{\text{-}\mathrm{gmod}}

\newcommand{\Irr}{{\rm Irr}}

\newcommand{\K}{\sfK}
\newcommand{\ex}{\mathrm{ex}}
\newcommand{\fr}{\mathrm{fr}}
\newcommand{\Kex}{{\K_\ex}}
\newcommand{\Kfr}{{\K_\fr}}

\newcommand{\seed}{\calS}

\newcommand{\g}{\mathfrak{g}}

\newcommand{\n}{\mathfrak{n}}

\newcommand{\Q}{\mathbb{Q}}
\newcommand{\Z}{\mathbb{Z}\ms{1mu}}

\newcommand{\al}{{\ms{1mu}\alpha}}

\newcommand{\la}{\lambda}
\newcommand{\be}{{\ms{1mu}\beta}}
\newcommand{\ga}{\gamma}
\newcommand{\ve}{\varepsilon}
\newcommand{\La}{\Lambda}

\newcommand{\rmL}{{\mathrm L}}
\newcommand{\rmG}{{\mathrm G}}

\newcommand{\wt}{{\rm wt}}
\newcommand{\PBW}{{\rm PBW}}

\newcommand{\de}{\mathfrak{d}}

\newcommand{\ii}{ \textbf{\textit{i}}}

\newcommand{\jj}{ \textbf{\textit{j}}}

\newcommand{\Hom}{\operatorname{Hom}}

\newcommand{\HOM}{\mathrm{H{\scriptstyle OM}}}

\newcommand{\tC}{\widetilde{C}}

\newcommand{\tB}{\widetilde{B}}

\newcommand{\tscrC}{\widetilde{\scrC}}

\newcommand{\tLa}{\widetilde{\Lambda}}

\newcommand{\uv}{{\underline{v}}}

\newcommand{\uL}{\underline{L}}

\newcommand{\frakS}{\mathfrak{S}}

\newcommand{\sfS}{\mathsf{S}}

\newcommand{\sfA}{\mathsf{A}}

\newcommand{\sfD}{\mathsf{D}}
\newcommand{\sfP}{\mathsf{P}}
\newcommand{\sfQ}{\mathsf{Q}}
\newcommand{\sfW}{\mathsf{W}}

\newcommand{\sfX}{\mathsf{X}}

\newcommand{\sfd}{\mathsf{d}}

\newcommand{\sfK}{\mathsf{K}}
\newcommand{\sfH}{\mathsf{H}}
\newcommand{\sfR}{\mathsf{R}}

\newcommand{\bbA}{\mathbb{A}}

\newcommand{\bbF}{\mathbb{F}}

\newcommand{\bfk}{\mathbf{k}}

\newcommand{\bfone}{\mathbf{1}}

\newcommand{\bfg}{\mathbf{g}}
\newcommand{\bfc}{\mathbf{c}}
\newcommand{\bfd}{\mathbf{d}}

\newcommand{\bfe}{\mathbf{e}}

\newcommand{\bfa}{\mathbf{a}}
\newcommand{\bfb}{\mathbf{b}}

\newcommand{\bfn}{\mathbf{n}}
\newcommand{\bfx}{\mathbf{x}}

\newcommand{\calD}{\mathcal{D}}
\newcommand{\calC}{\mathcal{C}}
\newcommand{\calQ}{\mathcal{Q}}

\newcommand{\calU}{\mathcal{U}}
\newcommand{\calA}{\mathcal{A}}
\newcommand{\calS}{\mathcal{S}}

\newcommand{\calK}{\mathcal{K}}

\newcommand{\calM}{\mathcal{M}}

\newcommand{\calT}{\mathcal{T}}

\newcommand{\calZ}{\mathcal{Z}}

\newcommand{\scrC}{\mathscr{C}}
\newcommand{\scrA}{\mathscr{A}}

\newcommand{\To}[1][{\hspace{2ex}}]{\xrightarrow{\,#1\,}}

\newlength{\mylength}
\newcommand*{\para}{%
  \rlap{\rotatebox{-30}{\rule[.05ex]{.4pt}{.77em}}}%
  \kern.04em%
  \rlap{\kern.36em\raisebox{0.649519052835em}{\rule{.6em}{.4pt}}}%
  \rule{.6em}{.4pt}\kern-.04em%
  \rotatebox{-30}{\rule[.05ex]{.4pt}{.77em}}}

\newcommand{\Rr}{\mathbf{r}}

\nc{\rmat}[1]{{\mathbf
r}_{\mspace{-2mu}\raisebox{-.5ex}{${\scriptstyle{#1}}$}}}

\newcommand{\wl}{\sfP}
\newcommand{\rl}{\sfQ}

\newcommand{\weyl}{\sfW}
\newcommand{\lan}{\langle}
\newcommand{\ran}{\rangle}

\newcommand{\isoto}[1][]{\mathop{\xrightarrow%
[{\raisebox{.3ex}[0ex][.3ex]{$\scriptstyle{#1}$}}]%
{{\raisebox{-.6ex}[0ex][-.6ex]{$\mspace{2mu}\sim\mspace{2mu}$}}}}}

\newcommand{\ee}{\end{enumerate}}
\newcommand{\bitem}{\begin{itemize}}
\newcommand{\eitem}{\end{itemize}}
\newcommand{\ben}{\begin{enumerate}[{\rm (1)}]}
\newcommand{\bnum}{\begin{enumerate}[{\rm (i)}]}
\newcommand{\bnump}{\begin{enumerate}[{\rm (i)$'$}]}
\newcommand{\bna}{\begin{enumerate}[{\rm (a)}]}
\newcommand{\bnA}{\begin{enumerate}[{\rm (A)}]}
\newcommand{\bc}{\begin{cases}}
\newcommand{\ec}{\end{cases}}

\newenvironment{myequation}
{\relax\setlength{\arraycolsep}{1pt}\begin{eqnarray}}
{\end{eqnarray}}
\newenvironment{myequationn}
{\relax\setlength{\arraycolsep}{1pt}\begin{eqnarray*}}
{\end{eqnarray*}}

\newcommand{\eq}{\begin{myequation}}
\newcommand{\eneq}{\end{myequation}}
\newcommand{\eqn}{\begin{myequationn}}
\newcommand{\eneqn}{\end{myequationn}}
\newcommand{\cl}{\colon}

\newcommand{\bl}{\bigl(}
\newcommand{\br}{\bigr)}
\newcommand{\qt}[1]{\quad\text{#1}}

\newcommand{\ba}{\begin{array}}
\newcommand{\ea}{\end{array}}

\nc{\on}{\operatorname}
\nc{\Aut}{\on{Aut}}
\nc{\wtl}{\wl}
\nc{\stt}[1]{\{#1\}}
\nc{\catC}{\scrC}
\nc{\one}{\bfone}
\nc{\ro}{{\rm(}}
\nc{\Proof}{\begin{proof}}
\nc{\QED}{\end{proof}}
\nc{\shs}{\mathcal{S}}
\nc{\GR}[1][{\calM}]{\rmG^R_{#1}}
\nc{\GL}[1][{\calM}]{\rmG^L_{#1}}

\usepackage[colorinlistoftodos]{todonotes}

\newcounter{myc}

\author[M. Kashiwara]{Masaki Kashiwara}
\thanks{The research of M.\ Kashiwara
was supported by Grant-in-Aid for Scientific Research (B) 20H01795,
Japan Society for the Promotion of Science.}
\address[M. Kashiwara]{%
Kyoto University Institute for Advanced Study, Research Institute
for Mathematical Sciences, Kyoto University, Kyoto 606-8502, Japan
\& Korea Institute for Advanced Study, Seoul 02455, Korea }
\email[M. Kashiwara]{masaki@kurims.kyoto-u.ac.jp}
\author[M. Kim]{Myungho Kim}
\address[M. Kim]{Department of Mathematics, Kyung Hee University, Seoul 02447, Korea}
\email[M. Kim]{mkim@khu.ac.kr}
\thanks{The research of M.\ Kim was supported by the National Research Foundation of
Korea (NRF) Grant funded by the Korea government(MSIP)
(NRF-2022R1F1A1076214 and NRF-2020R1A5A1016126).}
\author[S.-j. Oh]{Se-jin Oh}
\thanks{ The research of S.-j.\ Oh was supported by the Ministry of Education of the Republic of Korea and the National Research Foundation of Korea  (NRF-2022R1A2C1004045).}
\address[S.-j. Oh]{Department of Mathematics, Sungkyunkwan University, Suwon 16419, South Korea}
\email[S.-j. Oh]{sejin092@gmail.com}
\author[E. Park]{Euiyong Park}
\thanks{The research of E.\ Park was supported by the National Research Foundation of Korea (NRF) Grant funded by the Korea Government(MSIP)(NRF-2020R1F1A1A01065992 and NRF-2020R1A5A1016126).}
\address[E. Park]{Department of Mathematics, University of Seoul, Seoul 02504, Korea}
\email[E. Park]{epark@uos.ac.kr}

\title[Laurent family]{Laurent family of simple modules \\ over quiver Hecke algebra}



\date{June 27, 2023}

\begin{document}

\begin{abstract}
We introduce the notions of quasi-Laurent and
Laurent families of simple modules over quiver Hecke algebras of arbitrary symmetrizable types.
We prove that such a family plays a similar role of a cluster in the quantum cluster algebra theory and
exhibits a quantum Laurent positivity phenomenon to the basis of the quantum unipotent coordinate ring $\calA_q(\n(w))$, coming from
the categorification.
Then we show that the families of simple modules categorifying GLS-clusters are Laurent families by using the PBW-decomposition vector of a simple module $X$ and categorical interpretation of (co-)degree of $[X]$. As applications of such $\Z$-vectors, we define several skew symmetric pairings on arbitrary  pairs of simple modules, and investigate the relationships among the pairings and $\Lambda$-invariants of $R$-matrices in the quiver Hecke algebra theory.
\end{abstract}

\setcounter{tocdepth}{3} 

\maketitle
\tableofcontents

\section{Introduction}
A cluster algebra and its non-commutative  version
 quantum cluster algebra,
are introduced by Berenstein-Fomin-Zelevinsky~\cite{FZ02,BZ05} in an attempt to provide an algebraic and combinatorial framework for investigating the upper global basis of the quantum group.

The quantum cluster algebra $\scrA_q$ is a non-commutative $\Z[q^{\pm 1/2}]$-subalgebra
in the skew field $\Q(q^{1/2})(X_i)_{i \in \sfK}$ generated by   the   cluster variables, which are obtained from the initial cluster $\{  X_{i} \}_{i \in
\sfK}$ via the sequences of procedures, called \emph{mutations}. Even though the mutation involves non-trivial fraction, 
$\scrA_q$ is still contained in $\Z[q^{\pm 1/2}][X_i^{\pm1}]_{i \in \sfK}$
 with amazing reductions of fractions which is referred to as the \emph{quantum Laurent phenomenon} (\cite{BZ05}). The famous conjecture, which is not completely proved yet at this moment, is the \emph{quantum Laurent positivity} conjecture which asserts that every cluster variable is an element in $\Z_{\ge 0}[q^{\pm 1/2}][\sfX_i^{\pm1}]_{i \in \sfK}$ for any cluster $\{\sfX_i\}_{i \in \sfK}$. Note that the conjecture is proved in \cite{D18} (see also \cite{LS15,GHKK}) when $\scrA_q$
is of skew-symmetric type and is widely open when it is of non skew-symmetric type.

The notion of \emph{monoidal categorification} of (quantum) cluster algebra
is introduced by Hernandez-Leclerc in \cite{HL10} (see also \cite{KKKO18}) as the categorical framework for proving the conjecture as follows: A monoidal category $\calC$ with an  autofunctor $q$ is a monoidal categorification of $\scrA_q$, if {\rm (a)} $\bbA \otimes_{\Z[q^{\pm 1}]} K(\calC)$ $(\bbA \seteq\Z[q^{\pm 1/2}])$ is isomorphic to $\scrA_q$ and {\rm (b)} the cluster monomials of $\scrA_q$ are of the classes of real simple objects of $\calC$.
 Once $\calC$
is  a monoidal categorification of $\scrA_q$, then the conjecture for $\scrA_q$
follows since it can be interpreted as
the existence of a Jordan-H\"older series of an object.
 In \cite{KKKO18}, it is proved that the category
$\scrC_w$ over \emph{symmetric} quiver Hecke algebra $\sfR$
is a monoidal categorification
of the quantum unipotent coordinate ring $\calA_\bbA(\n(w))$ associated with
an element $w$ of the Weyl group $\weyl$
by using the $\Z$-invariant $\Lambda(M,N)$
of a pair of simple objects $M,N \in \scrC_w$.

For \emph{non-symmetric cases}, the monoidal categorification is still out of reach.
We know that $\scrC_w$ categorifies $\calA_\bbA(\n(w))$ as an algebra  (\cite{KL1,KL2,R08,Kimura12})
and $\calA_\bbA(\n(w))$ has a quantum cluster algebra structure (\cite{GLS13,GY17}) in \emph{every symmetrizable case}.
The quantum cluster algebra structure is skew-symmetric if the corresponding
generalized Cartan matrix is symmetric.
However,  we cannot prove that $\scrC_w$ is a
 monoidal categorification of $\calA_\bbA(\n(w))$ in \emph{non-symmetric cases}  due to the
  obstacle   that we do not know whether every simple module $M \in \scrC_w$ admits an affinization (\cite{KP18}) or not.  Note that the existence of affinizations guarantees to define R-matrices and the $\Z$-invariant $\Lambda(M,N)$.

\smallskip

In this paper, we study the quantum Laurent positivity for $\calA_\bbA(\n(w))$ of not necessarily symmetric type in the view point of the monoidal categorification.
More precisely, we show that  the basis of $\calA_\bbA(\n(w))$ corresponding to the simple modules in $\scrC_w$ exhibits  a quantum Laurent positivity phenomenon with respect to any \emph{quasi-Laurent family}, which is a central notion we introduce in this paper and plays the similar role of a cluster in the quantum cluster algebra theory.

The \emph{quasi-Laurent family} (resp.\ \emph{Laurent family}) $\calM=\{ M_j \}_{j \in J}$
consists of mutually commuting affreal simple modules in $\scrC_w$ satisfying additional conditions
(Definition~\ref{def: quasi-Laurent}).
Among others, the most important condition is that if a simple module $X$ commutes with all $M_j$, then there are monomials (i.e., convolution products) $\calM(\bold a)$ and $\calM(\bold b)$  in $\{M_j\}_{j\in J}$ such that $X \conv \calM(\bold a)$ is isomorphic to $\calM(\bold b)$.
We call the family  $\calM$ is \emph{Laurent} if  $\calM$ is \emph{maximal} in the sense that, if a simple module $X$ commutes with all $M_j$, then $X$
is isomorphic to a monomial $\calM(\bold b)$ in $\{M_j\}_{j\in J}$.

The main results of this paper are the followings:
\bnA
\item \label{it: Laurent_ph}
We show that if $\calM$ is a quasi-Laurent family in $\catC_w$, then the class $[X]$ in the Grothendieck ring $K(\catC_w)$ of any  simple object $X$ in $\catC_w$ can be written as a Laurent polynomial in $\{ [M_j] \}_{j \in J}$
whose coefficients belong to $\Z_{\ge0}[q^{\pm 1}]$ (Proposition~\ref{prop:Layrent}).
\item If $\calM$ is a monoidal seed in $\catC_w$, then
$\calM$ is a Laurent family.
\item \label{it: positivity_GLS}
In particular,
for any reduced sequence $\ii$ of $w$, we show that the family
$\calM^\ii \seteq \{ M(w_{\le k}^\ii \varpi_{i_k},\varpi_{i_k}) \}$ is a Laurent family and hence any class $[X]$ of a module $X$ in $\catC_w$ can be written as a Laurent polynomial in the unipotent quantum minors $D(w_{\le k}^\ii \varpi_{i_k},\varpi_{i_k})$ with coefficients in $\Z_{\ge0}[q^{\pm 1}]$. Note that  $D(w_{\le k}^\ii \varpi_{i_k},\varpi_{i_k})=[M(w_{\le k}^\ii \varpi_{i_k},\varpi_{i_k}) ]$ and we call $\{D(w_{\le k}^\ii \varpi_{i_k},\varpi_{i_k})\}$ the \emph{GLS-seed associated with $\ii$} (Proposition~\ref{prop: Laurent commuting}).

\item We show that if $\calM$ is a quasi-Laurent family, then the class $[X]$ of a simple module $X$ in $\catC_w$ is pointed and copointed with respect to the partial order $\preceq_\calM$. That is, the set of degrees of the monomials appearing in the  Laurent expansion of  $[X]$ with respect to $\calM$ has a unique maximal element and a unique minimal element with respect to $\preceq_\calM$.
We define vectors $\bfg_\calM^R(X)$ and $\bfg_\calM^L(X) \in \Z^{\oplus J}$  as the maximal and the minimal element, respectively.

\item
Each  quasi-Laurent family $\calM$ also induces new $\Z$-values
$\rmG_\calM^R(X,Y)$ and
  $\rmG_\calM^L(X,Y)$ for \emph{any} pair of simple modules $X$ and $Y$ which coincides with $\La(X,Y)$ provided $X,Y$ commutes and one of them is affreal.
\ee

To our best knowledge, the positivity result in ~\eqref{it: positivity_GLS} is new. 
We can understand ~\eqref{it: Laurent_ph} that a quasi-Laurent family is a generalization of a cluster in the categorical view point, and
that the positivity conjecture  can be extended to elements corresponding to simple modules in all skew-symmetrizable types.

\smallskip

In \cite{FZ07,Qin17},
Fomin-Zelevinsky and Qin defined a pointed (resp. co-pointed) element $\bfx$ in a cluster algebra and its degree $\mathbf{deg}_{\calS}(\bfx) \in \Z^{\oplus \sfK}$ (resp.\ codegree $\mathbf{codeg}_{\calS}(\bfx) \in \Z^{\oplus \sfK}$) depending on the choice of a seed $\calS$ (see also \cite{Qin20} for codegree and \cite{Tran} for degree elements in a quantum cluster algebra).
With a fixed choice of a seed,
it is proved in \cite{Tran} that every cluster monomial is pointed, and in \cite{DWZ10,GHKK}
that cluster monomials are determined by their degrees.

For a given quasi-Laurent family $\calM$ and a simple module $X \in \scrC_w$,
we define vectors $\bfg_\calM^R(X),\bfg_\calM^L(X) \in \Z^{\oplus J }$
in Definition~\ref{def: g-vectors} by using the $\Z^{\oplus J}_{\ge 0}$-vectors
in Lemma~\ref{lem: a,a',b,b'} and guaranteeing its well-definedness in Lemma \ref{lem:well-defined}.
We then prove that,
for every simple module $X \in \scrC_w$, the element $[X]$ in $\calA_\bbA(\n(w))$is  is (co-)pointed with respect to the GLS-seed $\calS^\ii$ and
that $\bfg_{\calM^{\ii}}^R(X)$ and $\bfg_{\calM^{\ii}}^L(X)$ coincide with $\mathbf{deg}_{\calS^\ii}([X])$ and $\mathbf{codeg}_{\calS^\ii}([X])$, respectively.

Utilizing the vectors $\bfg_\calM^R(X)$ and $\bfg_\calM^L(X)$, we define skew-symmetric $\Z$-valued forms $\rmG_\calM^R(-,-)$
and $\rmG_\calM^L(-,-)$ on the pairs $(X,Y)$ of simple modules.  Then we compare
$\rmG_\calM^R(X,Y)$ and $\rmG_\calM^L(X,Y)$ with the $\Z$-invariant $\La(X,Y)$

when the pair of simple module $(X,Y)$ admits the $\Z$-invariant $\La(X,Y)$.
It is proved in Proposition~\ref{prop: lower bound} that $\rmG_\calM^R(X,Y)$ and $\rmG_\calM^L(X,Y)$ give lower bounds of $\La(X,Y)$,
and in Proposition~\ref{prop: Gms Lambda} that
$\rmG_\calM^R(X,Y)=\rmG_\calM^L(X,Y)=\La(X,Y)$ when $(X,Y)$ is a commuting pair. Here we would like to emphasize that {\rm (1)} $\rmG_\calM^R(X,Y)$ and $\rmG_\calM^L(X,Y)$
are defined even for pairs $(X,Y)$ we do not know whether they admit $\La(X,Y)$ or not, and {\rm (2)} the $\Z$-values $\rmG_\calM^R(X,Y)$ and $\rmG_\calM^L(X,Y)$ \emph{do depend} on the choice of $\calM$ as (co-)degree does on the one of seeds (Remark~\ref{rmk: bfg depend}).

\smallskip

This paper is organized as follows. In Section~\ref{sec: Preliminaries}, we give preliminaries.
In Section~\ref{sec: Laurent family}, we define the notions of quasi-Laurent and Laurent families, and investigate their properties. Then we define $\bfg_{\calM}^R(X)$ and $\bfg_{\calM}^L(X)$, and prove that $\bfg_{\calM}^R(X)$ and $\bfg_{\calM}^L(X)$ determine the isomorphism class of $X$. In Section~\ref{sec: PBW and GLS}, we prove that $\calM^\ii$ is Laurent by studying PBW-decomposition vectors
 of simple modules. In Section~\ref{sec: g-vector}, we define the skew symmetric pairings on  pairs of simple modules  and investigate the relationships among the pairings and $\Lambda$-invariants.

\medskip \noindent
{\bf Acknowledgments.}\
The second, third and fourth authors gratefully acknowledge for the hospitality of RIMS (Kyoto University) during their visit in 2023.

\begin{convt}  Throughout this paper, we use the following convention.
\bnum
\item For a statement $\mathtt{P}$, we set $\delta(\mathtt{P})$ to be $1$ or $0$
depending on whether $\mathtt{P}$ is true or not.
As a special case, we use the notation $\delta_{i,j} \seteq \delta(i=j)$ (Kronecker's delta).
\item For integers $a,b \in \Z$, we set
\begin{align*}
&[a,b] \seteq \{ x\in \Z \; | \;  a \le x \le b\}
\end{align*}
We refer to the subset as an \emph{interval} and understand an empty set if $a>b$.
\item Let $\bfx=(x_j )_{ \ j \in J }$ be a family parameterized by an index set $J$. Then for any $j \in J$, we set
$$    (\bfx)_j \seteq x_j.$$
\ee
\end{convt}

\section{Preliminaries} \label{sec: Preliminaries}

In this preliminary section, we briefly review the basic materials of this paper.
We refer the reader to \cite{BZ05,FZ07,KL1,R08,Kimura12,GLS13,KKKO18,GY17,KO17,KP18,KKOP18,KK19,GHKK} for more details.

\subsection{Quantum cluster algebras} \label{sec: quantum cluster algebra} Fix a finite index set $\K = \Kex \sqcup \Kfr$ with a decomposition into the set $\Kex$ of exchangeable indices and the set $\Kfr$ of frozen indices.
Let $L = (l_{ij})_{i,j \in \K}$ be a skew-symmetric integer-valued matrix and let $q$ be an indeterminate. We set $\bbA \seteq \Z[q^{\pm 1/2}]$
where $q^{1/2}$ denotes the formal square root of $q$.

\begin{definition} We define the \emph{quantum torus} $\calT(L)$ to be the $\bbA$-algebra generated by a  finite  family of elements $\{X^{\pm 1}_k \}_{k \in\K}$ subject to the following defining relations:
$$
X_jX_j^{-1} = X_j^{-1}X_j=1 \quad\text{and}\quad X_iX_j = q^{l_{ij}}X_jX_i \quad \text{ for $i,j \in \K$.}
$$
\end{definition}

For $\bfa = (\bfa_i)_{i \in \K} \in \Z^\K$, we define the element $X^\bfa$ of $\calT(L)$ as
$$
X^\bfa = q^{\frac{1}{2} \sum_{i>j}\bfa_i\bfa_j l_{ij}} \prod^{\Lto}_{i \in \K} X_i^{\bfa_i}.
$$
Here  $\overset{\Lto}{\prod}_{i \in \K} X_i^{\bfa_i} \seteq X_{i_1}^{\bfa_{i_1}} \cdots  X_{i_r}^{\bfa_{i_r}}$,
where $\K=\{ i_1,\ldots,i_r \}$  with a total order $i_1 < \cdots <i_r$. Note that $X^\bfa$ does not depend on the choice of a total order $<$ on $\K$. Then $\{ X^\bfa \mid\bfa \in \Z^{\K} \}$ forms an $\bbA$-basis of $\calT(L)$. Since
$\calT(L)$ is an Ore domain, it is embedded into the skew field of fractions $\bbF(\calT(L))$.

Let $\tB = (b_{ij})_{i \in \K,j\in \Kex}$ be an integer-valued $\K \times \Kex$-matrix whose principal part $B = (b_{ij})_{i,j\in \Kex}$ is skew-symmetrizable; i.e., there exists a diagonal matrix $D$ with a positive integer entries
such that $DB$ is skew-symmetric. Such a matrix $\tB$ is called an \emph{exchange matrix}.  We say that a pair $(L,B)$ is \emph{compatible} if
$$  \sum_{k \in \K} b_{ki} l_{kj}  = d_i\delta_{i,j}\qt{for any $i\in \Kex$ and $j\in\K$} $$
for some positive integers $\{ d_i \}_{i \in \Kex}$.  We call the triple $\calS = (\{ X_k \}_{k \in \K}, L , \tB )$ a \emph{quantum seed} in the quantum torus $\calT(L)$ and $\{ X_k \}_{k \in \K}$ a \emph{quantum cluster}. 

For $k \in \Kex$, the \emph{mutation} $\mu_k(L,\tB) \seteq (\mu_k(L),\mu_k(B))$ of a compatible pair $(L,\tB)$ \emph{in a direction $k$} is defined in a combinatorial way (see~\cite{BZ05}).
Note that (i) the pair $(\mu_k(L),\mu_k(B))$ is also compatible with the same positive integers $\{d_i\}_{i \in \K}$ and (ii) the operation $\mu_k$ is an involution; i.e.,
$ \mu_k(\mu_k(L,\tB)) = (L,\tB)$.
We define an isomorphism of $\Q(q^{1/2})$-algebras $\mu_k^* \colon \bbF(\calT(\mu_kL)) \isoto  \bbF(\calT(L))$ by
$$
\mu_k^*(X_j) \seteq \bc
X^{\bfa'} + X^{\bfa''} & \text{ if } j =k, \\
X_j  & \text{ if } j \ne k,
\ec
$$
where
$$
\bfa'_i = \bc
-1 & \text{ if } i =k, \\
\max(0,b_{ik}) & \text{ if } i \ne k,
\ec
\text{ and } \ \
\bfa''_i = \bc
-1 & \text{ if } i =k, \\
\max(0,-b_{ik}) & \text{ if } i \ne k.
\ec
$$
Then the \emph{mutation $\mu_k(\seed)$ of the quantum seed $\seed$ in a direction $k$} is defined to be the triple $( \{ X_i\}_{i \ne k} \sqcup \{ \mu_k^*(X_k) \}, \mu_k(L), \mu_k(\tB))$.

For a quantum seed $\seed=( \{X_k\}_{k\in \K},L,\tB)$, an element in $\bbF(\calT(L))$ is called a \emph{quantum cluster variable} (resp.\ \emph{quantum cluster monomial}) if it is of the form
$$
\mu_{k_1}^* \cdots \mu_{k_\ell}^*(X_j), \quad \text{$($resp.\ } \mu_{k_1}^* \cdots \mu_{k_\ell}^*(X^\bfa) )
$$
for some finite sequence $(k_1,\ldots,k_\ell) \in \K_\ex^\ell$ $(\ell \in \Z_{\ge 0})$ and $j \in \K$ (resp.\ $\bfa \in \Z_{\ge 0}^\K$). For a quantum seed $\calS = (\{ X_k \}_{k \in \K}, L,\tB )$, the \emph{quantum cluster algebra}
$\scrA_q(\calS)$ is the $\bbA$-subalgebra of $\bbF(\calT(L))$ generated by all the quantum cluster variables. Note that   $\scrA_q(\calS)\simeq\scrA_q(\boldsymbol{\mu}(\calS))$ for any sequence $\boldsymbol{\mu}$ of mutations.

\smallskip

The \emph{quantum Laurent phenomenon}, proved by Berenstein-Zelevinsky in \cite{BZ05}, tells that the quantum cluster algebra $\scrA_q(\calS)$ is contained in $\calT(L)$ indeed.

\smallskip

For a quantum seed $\seed$ with a compatible pair $(L,\tB)$, an element $\bfx \in \calT(L)$ is called \emph{pointed} (resp.\  \emph{copointed}) if it is of the following form:
\begin{align} \label{eq: def degree and codegree}
\bfx =  q^{a} X^{\bfg^R} + \sum_{ \bfc \in \Z_{\ge 0}^{\Kex} \setminus \{ 0 \} } p_\bfc X^{\bfg^R+\tB \bfc} \ \ \text{$($resp.\  } \bfx =  q^{a} X^{\bfg^L} + \sum_{ \bfc \in \Z_{\le 0}^{\Kex} \setminus \{ 0 \} } p_\bfc X^{\bfg^L+\tB \bfc})
\end{align}
for some $a \in \dfrac{1}{2}\Z$, $\bfg^R \in \Z^\K$ (resp.\  $\bfg^L \in \Z^\K$) and $p_\bfc \in \bbA$. In this case, we call $\bfg^R$ the \emph{degree} (resp.\  \emph{codegree}) of the pointed (resp.\  copointed)
element $\bfx$ and denote it by $\mathbf{deg}_{\calS}(\bfx)$ (resp.\  $\mathbf{codeg}_{ \calS}(\bfx)$). The degree (resp.\  codegree) is often the called \emph{$g$-vector} (resp.\  \emph{dual $g$-vector}) of $\bfx$  (see \cite[Definition 3.1.4]{Qin17} and \cite[Definition 3.1.3]{Qin20}).
It is worth to remark that the notion of $g$-vector (resp.\  dual $g$-vector) \emph{does depend on} the compatible pair $(L,\tB)$ and hence   on  the seed $\seed$. It is proved in \cite[Theorem 5.3]{Tran} that every quantum cluster monomial in $\scrA_q(\seed)$ is pointed.

We say that an  $\bbA$-algebra $R$ has a
{\em quantum cluster algebra structure}
if there exists a quantum seed $\seed$ and an $\bbA$-algebra isomorphism $\Omega: \scrA_q(\seed) \isoto R$.
In the case, a {\em quantum seed of $R$} refers to the image of a quantum seed in $\scrA_q(\seed)$, which is obtained by a sequence of mutations.

\subsection{Quantum unipotent coordinate rings}  \label{subsec: Quantum unipotent coordinate ring}

Let $I$ be an index set. A \emph{Cartan datum} $(\sfA,\wl,\Pi,\wl^\vee,\Pi^\vee)$ consists of
\bna
\item a symmetrizable Cartan matrix $\sfA=(a_{i,j})_{i,j \in I}$, i.e., $\sfD \sfA$ is symmetric for a diagonal matrix $\sfD = {\rm diag}(\sfd_i \mid i \in I)$
with $\sfd_i\in\Z_{>0}$,
\item a free abelian group $\wl$, called the \emph{weight lattice},
\item $\Pi = \{\al_i \mid i \in I \} \subset \wl$ , called the set of \emph{simple roots},
\item $\Pi^\vee = \{ h_i \mid i \in I \} \subset \wl^\vee \seteq \Hom(\wl,\Z)$ , called the set of \emph{simple coroots},
\item a $\Q$-valued symmetric bilinear form $( \cdot,\cdot)$  on $\wl$,
\ee
satisfying the standard properties (see~\cite[\S 1.1]{KKKO18} for instance).
Here we take $\sfD ={\rm diag}(\sfd_i \mid i \in I)$ such that  $\sfd_i \seteq(\al_i,\al_i)/2\in\Z_{>0}$  $(i\in I)$ in this paper.

For $i \in I$, we choose $\varpi_i \in \wl$ such that $\lan h_i,\varpi_j \ran =\delta_{ij}$ for any $j \in I$ and call it the \emph{$i$-th fundamental weight.}
The free abelian group $\rl \seteq \bigoplus_{i \in I} \Z \al_i$ is called the \emph{root lattice} and we set $\rl^+ = \sum_{i \in I}\Z_{\ge 0} \al_i \subset \rl$ and $\rl^- = \sum_{i \in I}\Z_{\le 0} \al_i \subset \rl$.
We denote by $\Delta$ the set of \emph{roots} and by $\Delta^\pm$ the set of \emph{positive} roots (resp.\ \emph{negative} roots).
For $\be \in \sum_{i \in I} m_i\al_i \in \rl^+$, we set $|\be| \seteq \sum_{i \in I} m_i$, ${\rm supp}(\be) \seteq\{  i \in I \mid  m_i \ne 0 \}$ and $I^\be \seteq \{ \nu = (\nu_1, \ldots,\nu_{|\be|}) \in I^{|\be|} \mid  \al_{\nu_1}+ \cdots + \al_{\nu_{|\be|}} =\be \}$.
Note that the symmetric group $\frakS_{|\be|} = \lan r_{1},\ldots, r_{|\be|}  \ran$ acts on $I^\be$ by the place permutations.

Let $\g$ be the Kac-Moody algebra associated with the Cartan datum $(\sfA,\wl,\Pi,\wl^\vee,\Pi^\vee)$, and $\weyl$
the \emph{Weyl group} of $\g$. It is generated by the simple reflections $s_i\in\Aut(\wtl)$ $(i \in I)$ defined by $s_i(\la) = \la - \lan h_i,\la \ran \al_i$ for
$\la \in \wl$. For a sequence $\ii=(i_1,\ldots,i_r) \in I^r$, we call it  a \emph{reduced sequence} of $w \in \weyl$ if $s_{i_1} \ldots s_{i_r}$ is a reduced expression of $w$. For $w,v \in \weyl$, we write $w \ge v$ if there is a reduced sequence of $v$ which appears in a reduced sequence of $w$ as a subsequence.

For $\la,\mu \in \wl$, we write $\la \preceq \mu$ if there exists a sequence of
real positive roots $\be_k$ $(1 \le k \le l)$ such that
$\la = s_{\be_l} \cdots s_{\be_1} \mu$ and  $(\be_k,s_{\be_{k-1}} \cdots s_{\be_1}\mu)>0$ for $1 \le k \le l$. When $\La \in\wl^+$ and $\la,\mu \in \weyl\La$
the relation $\la \preceq \mu$ holds if and only if there exist $w,v \in \weyl$
such that $\la = w\La$, $\mu = v\La$ and $v \le w$.

Let $\calU_q(\g)$ be the quantum group of $\g$ over $\Q(q^{1/2})$, generated by $e_i,f_i$ $(i \in I)$ and $q^h$ $(h \in \wl^\vee)$.
We denote by $\calU_q^+(\g)$ the subalgebra of $\calU_q(\g)$ generated by $e_i$ and $\calU_{\bbA}^+(\g)$ the $\bbA$-subalgebra of
$\calU_q(\g)^+$ generated by $e_i/[n]_i !$ $(i \in I, \; n \in \Z_{>0})$, where
$$
q_i \seteq q^{\sfd_i}, \quad [k]_i = \dfrac{q_i^k-q_i^{-k}}{q_i - q_i^{-1}} \quad\text{and}\quad [k]_i! = \prod_{s=1}^k [s]_i.
$$
Set
$$ \calA_q(\n) = \soplus_{\beta \in  \rl^-} \calA_q(\n)_\beta \quad \text{ where } \calA_q(\n)_\beta \seteq \Hom_{\Q(q^{1/2})}(\calU^+_q(\g)_{-\beta}, \Q(q^{1/2})),$$
where $\calU^+_q(\g)_{-\beta}$ denote the $(-\be)$-root space of $\calU^+_q(\g)$.
Then  $\calA_q(\n)$ also has an algebra structure and is called  the \emph{quantum unipotent coordinate ring} of $\g$. We denote by $\calA_{\bbA}(\n)$ the $\bbA$-submodule of
$\calA_q(\n)$ generated by $\uppsi \in \calA_q(\n)$ such that $\uppsi(\calU_{\bbA}^+(\g)) \subset \bbA$.
Then, $\calA_q(\n)$ is an $\bbA$-subalgebra with a $\calU_{\bbA}^+(\g)$-bimodule structure.

For each $\la \in \wl^+ \seteq \sum_{i \in I} \Z_{\ge 0} \varpi_i$ and Weyl group elements $w,w' \in \weyl$, we can define a specific homogeneous element $D(w\la,w'\la)$ of $\calA_\bbA(\n)$, called a \emph{unipotent quantum minor}
(see for example \cite[Section 9]{KKKO18}).

For $w \in \weyl$, we denote by $\calA_\bbA(\n(w))$ the $\bbA$-submodule of $\calA_{\bbA}(\n)$ consisting of elements $\uppsi$ such that
$e_{i_1} \cdots e_{i_{|\be|}}\uppsi =0$ for any $\be \in  \rl^+ \setminus w \rl^-$ and $(\nu_{i_1},\ldots,\nu_{i_{|\be|}}) \in I^\be$. Then it is an $\bbA$-subalgebra and
we call it the
 \emph{quantum unipotent coordinate ring associated with $w$}.

For a reduced sequence $\ii=(i_1,\ldots,i_r)$ of $w \in \weyl$ and $1 \le k \le r$, define $w^{\ii}_{\le k}=s_{i_1} \cdots  s_{i_k}$ and $w^\ii_{< k}=s_{i_1}\cdots s_{i_{k-1}}$. Then $\calA_\bbA(\n(w))$ is generated by the set
of unipotent quantum minors $\{ D(w^\ii_{\le k}\varpi_{i_{k}}, w^\ii_{< k}\varpi_{i_{k}}) \mid 1 \le k \le r \}$ as an $\bbA$-algebra.

It is proved in \cite{GLS13,GY17,KKKO18,Qin20}
that $\calA_\bbA(\n(w))$ has a quantum cluster algebra structure, one of whose quantum seed  $\seed^\ii$
can be obtained from a reduced sequence $\ii=(i_1,\ldots,i_r)$ of $w$. To introduce   $\seed^\ii$, we need preparations:

\smallskip

Let $\jj=(j_1,\ldots,j_l)$ be a sequence of indices in $I$.
For $1 \le k \le l$  and $j \in I$, we set
\eqn
k^\jj_+&&\seteq \min(\stt{u\mid k<u \le l,\; j_u=j_k} \cup \{l+1\}), \\
k^\jj_-&&\seteq \max(\stt{u\mid 1\le u < k,\;j_u=j_k} \cup \{0\}).
\eneqn
We also set
$$
k^\jj_{\min} \seteq  \min\stt{u\mid 1\le u\le k ,\ j_u=j_k}  \quad\text{and}\quad k^\jj_{\max} \seteq  \max\stt{u\mid  k\le u\le l,\ j_u=j_k}.
$$
We sometimes drop $^\jj$ in the above notations if there is no danger of confusion.

Take $\K=[1,r]$ as an index set and  decompose $\K$ into
$$ \Kfr = \{ k \mid  1\le  k \le r, \; k^\ii_+ = r+1 \} \quad\text{ and }\quad \Kex \seteq  \K \setminus \Kfr.$$
We define the $\Z$-valued $\K\times \Kex$ matrix $\tB^\ii =(b^\ii_{st})_{s\in\K,t\in \Kex}$ and the $\Z$-valued skew-symmetric $\K\times \K$ matrix $L^\ii=(l^\ii_{st})_{s,t \in \K}$ as follows:
\begin{equation}\label{eq: tBi Li}
\begin{aligned}
b^\ii_{st} &=
\bc
\pm 1 & \text{if $s=t^\ii_\pm$,} \\
-a_{i_s,i_t}& \text{if  $s<t<s^\ii_+<t^\ii_+$,}\\
a_{i_s,i_t}&\text{if  $t<s<t^\ii_+<s^\ii_+$,} \\
0 & \text{otherwise},
\ec
  \\
l^\ii_{st} & = (\varpi_{i_s}-w^\ii_{\le s}\varpi_{i_s}, \varpi_{i_t}+w^\ii_{\le t}\varpi_{i_t}  ) \ \  \text{ for } s<t.
\end{aligned}
\end{equation}

Then the quantum seed $\seed^\ii$ of $\calA_\bbA(\n(w))$ is given as follows:
\begin{align}\label{eq: seed ii}
\seed^\ii \seteq \Bigl( \stt{ q^{c^\ii_k} D(w^\ii_{\le k}\varpi_{i_k},\varpi_{i_k}  )
}_{k \in \K}, L^\ii, \tB^\ii \Bigr),
\end{align}
where
$c^\ii_s = \frac{1}{4} (\varpi_{i_s}-w^\ii_{\le s}\varpi_{i_s} ,\varpi_{i_s}-w^\ii_{\le s}\varpi_{i_s} ) \in \Z/2$.
Note that
$(L^\ii\tB^\ii)_{ab} = -2 \sfd_{i_a} \times \delta_{a,b}$ for $(a,b)\in\K\times \Kex$,
$\wt(D(w^\ii_{\le k}\varpi_{i_k},\varpi_{i_k}  )) = - \varpi_{i_s} + w^\ii_{\le s}\varpi_{i_s}$, and $$\big\{ q^{c^\ii_k} D(w^\ii_{\le k}\varpi_{i_k},\varpi_{i_k}) = q^{c^\ii_k} D(w\varpi_{i_k},\varpi_{i_k} ) \mid k \in \Kfr \big\}$$ forms the set of frozen variables of the quantum cluster algebra
$\calA_\bbA(\n(w))$.
We call $\seed^\ii$ the {\em GLS-seed} (associated with $\ii$).
\smallskip

We set $\calD(w) \seteq  \{ q^m D(w\varpi,\varpi) \mid m \in \Z/2, \; \varpi \in \wl^+ \}$.
Then it is well-known that $\calD(w)$ consists of $q$-central elements of  $\calA_\bbA(\n(w))$  and hence forms an Ore set. We denote by $\calA_\bbA(\n^w)$
the quotient ring of
$\calA_\bbA(\n(w))$  by the Ore set $\calD(w)$.
Then $\calA_\bbA(\n^w)$ has also the quantum cluster algebra structure
with the \emph{invertible} frozen variables $\{ q^{c^\ii_k} D(w^\ii_{\le k}\varpi_{i_k},\varpi_{i_k})  \}_{k \in \Kfr}$
in the sense of \cite{BZ05}.

\subsection{Quiver Hecke algebras and categorifications} \label{subsec: Quiver Hecke algebras}

Let $\bfk$ be a base field. For $i,j \in I$, we choose polynomials $\calQ_{i,j}(u,v) \in \bfk[u,v]$ such that (a) $\calQ_{i,j}(u,v) =\calQ_{j,i}(v,u)$
and (b) each $\calQ_{i,j}(u,v)$ is of the following form:
$$
\calQ_{i,j}(u,v) =  \delta(i\ne j) \hspace{-6ex} \sum_{p(\al_i,\al_i)+q(\al_j,\al_j)=-2(\al_i,\al_j)}  \hspace{-6ex} t_{i,j;p,q} u^p v^q \quad \text{ where } \ \ t_{i,j;-a_{i,j},0} \in \bfk^\times.
$$

For a Cartan datum $(\sfA,\wl,\Pi,\wl^\vee,\Pi^\vee)$ and $\be \in \rl^+$, the \emph{quiver Hecke algebra $R(\be)$  associated with } $(\calQ_{i,j})_{i,j \in I}$ is the $\Z$-graded algebra over $\bfk$
generated by the elements
$$\{ e(\nu) \}_{\nu \in I^\be} ,  \ \   \{ x_k \}_{1 \le k \le |\be|},  \ \  \{ \tau_m \}_{1 \le m <|\be|}$$ subject to certain defining relations (see~\cite[Definition 1.1]{KKOP21} for instance).
Note that the
$\Z$-grading of $R(\be)$ is determined by the degrees of following elements:
$$
\deg(e(\nu))=0, \quad \deg(x_ke(\nu)) = (\al_{\nu_k},\al_{\nu_k}), \quad\text{and}\quad \deg(\tau_me(\nu)) = -(\al_{\nu_m},\al_{\nu_{m+1}}).
$$
We say that $R(\be)$ is \emph{symmetric} if $\calQ_{i,j}(u,v) \in \bfk[u-v]$ for $i,j \in {\rm supp}(\be)$.

\smallskip

We denote by $R(\be)\gmod$ the category of finite-dimensional graded $R(\be)$-modules with homomorphisms of degree $0$. For $M \in R(\be)\gmod$, we set $\wt(M) \seteq -\be \in \rl^-$.  Note that there exists
the \emph{degree shift functor}, denoted by $q$, such that $(qM)_n = M_{n-1}$ for $M=\soplus_{k \in \Z} M_k \in R(\be)\gmod$.

Throughout this paper, we usually deal with graded $R(\be)$-modules $(\be \in \rl^+)$ and sometimes skip grading shifts. Thus we usually say that $M$ is an $R$-module instead of saying that $M$ is
a graded $R(\be)$-module and $f\cl M \to N$ is a homomorphism if $f\cl q^aM \to N$ is a morphism in $R(\be)\gmod$. We set
$$
\HOM_{R(\be)}(M,N) \seteq \soplus_{a \in \Z} \HOM_{R(\be)}(M,N)_a
$$
with $\HOM_{R(\be)}(M,N)_a\seteq\Hom_{R(\be)\gmod}(q^aM,N)_a$.
We write $\deg(f) \seteq a$ for an $f \in \HOM_{R(\be)}(M,N)_a$.

For an $R(\be)$-module $M$ and an $R(\ga)$-module $N$, we can obtain $R(\be+\ga)$-module
$$ M \conv N \seteq R(\be+\ga) e(\be,\ga) \tens_{R(\be) \otimes R(\ga)} (M \otimes N).$$
where $e(\be,\ga) \seteq \sum_{\nu \in I^\be, \nu' \in I^\ga} e(\nu * \nu') \in R(\be+\ga)$. Here
$\nu*\nu'$ denotes the concatenation of $\nu$ and $\nu'$, and $\conv$ is called the  \emph{convolution product}.
We say that two simple $R$-modules $M$ and $N$ \emph{strongly commute} if $M \conv N$ is simple. If a simple module $M$ strongly commutes with itself, then $M$ is called \emph{real}. A simple $R$-module $M$ is said to be \emph{prime} if there are no non-trivial simple $R$-modules $N_1$ and $N_2$ such that $M \simeq N_1 \conv N_2$.

Set $R\gmod \seteq \soplus_{\be \in \rl^+}R(\be)\gmod$. Then the category $R\gmod$ is a monoidal category with the tensor product $\conv$ and the unit object $\mathbf{1} \seteq \bfk \in R(0)\gmod$. Hence the Grothendieck group
$K(R\gmod )$ has the $\Z[q^{\pm 1}]$-algebra structure derived from $\conv$ and the degree shift functors $q^{\pm1}$.

For a monoidal abelian subcategory $\calC$ of $R\gmod$ stable by grading shifts,
we set
$$ \calK_\bbA(\calC) \seteq \bbA \tens_{\Z[q^{\pm1}]} K(\calC ),$$
where $K(\calC )$ denotes the Grothendieck ring of $\calC$. For a subcategory $\scrC$ of
$R\gmod$, we denote by $\Irr(\scrC)$ the set of the isomorphism classes of \emph{self-dual} modules in $\scrC$.
Note that $\Irr(R\gmod)$  forms an $\bbA$-basis of  $\calK_\bbA(R\gmod)$.

It is proved in \cite{KL1,KL2,R08}
that there exists an $\bbA$-algebra isomorphism
\begin{align}\label{eq: categorification}
 \Upomega \cl    \calK_{\bbA}(R\gmod)   \isoto \calA_{\bbA}(\n).
\end{align}

\begin{proposition}[{\cite[Proposition 4.1]{KKOP18}}]  \label{prop: determinnatial}
For $\varpi \in \wl^+$ and $\mu,\zeta \in \weyl \varpi$ with $\mu \preceq  \zeta$, there exists a self-dual real simple $R( \zeta-\mu)$-module $M(\mu,\zeta)$ such that
$$\Upomega ([M(\mu,\zeta)]) = D(\mu,\zeta).$$
Here, $[M(\mu,\zeta)]$ denotes the isomorphism class of $M(\mu,\zeta)$ which is
called the \emph{determinantial module} associated with $D(\mu,\zeta)$.
\end{proposition}

For an $R(\be)$-module $M$, we define
\begin{align*}
W(M) & \seteq \{ \ga \in \rl^+ \cap (\be-\rl^+) \mid e(\ga,\be-\ga)M \ne 0 \}, \\
W^*(M) & \seteq \{ \ga \in \rl^+ \cap (\be-\rl^+) \mid e(\be-\ga,\ga)M \ne 0 \}.
\end{align*}

An ordered pair $(M, N)$ of $R$-modules is called \emph{unmixed} (\cite[Definition 2.5]{TW16}) if
$$W^*(M) \cap  W(N) \subset \{0\}.$$

For $w \in \weyl$, we denote by $\scrC_w$ the full subcategory of $R\gmod$ whose objects $M$ satisfy $W(M) \subset \rl^+ \cap w\rl^-$. Then the category $\scrC_w$ is the smallest monoidal abelian category of $R\gmod$ which (i) is stable under taking subquotients, extensions, grading shifts
and  (ii)  contains $\left\{ S^\ii_k \seteq M(w^\ii_{\le k}\varpi_{i_k},w^\ii_{< k}\varpi_{i_k})  \mid 1 \le k \le r \right\}$ for any reduced sequence $\ii$ of $w$. We call
$S^\ii_k$
the $k$-th \emph{cuspidal module associated with $\ii$}. Defining $\be_k^\ii \seteq w^\ii_{<k} \al_{i_k}$ for $1 \le k \le r$, one can see that
$\stt{\be_k^\ii \mid 1\le k\le r}=\Delta^+ \cap w \Delta^-$, and $-\wt(S^\ii_k) = \be^\ii_k$.
Then we have (\cite[Section 4]{Kimura12})
$$
 \Upomega( \calK_\bbA(\scrC_w ) ) \simeq \calA_\bbA(\n(w)).
$$

\subsection{$R$-matrices and affreal simple modules} For $\be \in \rl^+$ and $i \in I$, let
\begin{align*} 
\mathfrak{p}_{i,\be} = \sum_{ \eta\,  \in I^\be}  \big( \prod_{a \in [1,|\be|]; \  \eta_a  = i} x_a \big) e(\eta) \in \calZ(R(\be)),
\end{align*}
where $\calZ(R(\be))$ denotes the center of $R(\be)$.

\begin{definition}[{\cite[Definition 2.2]{KP18}}]
For an $R(\be)$-module $M$, we say that $M$ \emph{admits an affinization} if  there exists an $R(\be)$-module $\widehat{M}$
satisfying the condition: there exists an endomorphism $z_{\widehat{M}}$ of degree $t  \in \Z_{>0}$ such that $\widehat{M}/z_{\widehat{M}}\widehat{M} \simeq M$ and
\bnum
\item $\widehat{M}$ is a finitely generated free modules over the polynomial ring $\bfk[z_{\widehat{M}}]$,
\item $\mathfrak{p}_{i,\be}  \widehat{M} \ne 0$ for all $i \in I$.
\ee
We say that a simple $R(\be)$-module $M$ is \emph{affreal} if $M$ is real and admits an affinization.
\end{definition}

It is known that any $M \in R(\be)\gmod$ admits an affinization if $R(\be)$ is symmetric. However,  when $R(\be)$ is not symmetric, it is widely open whether an $R(\be)$-module  $M$ admits an affinization or not.

\begin{theorem} [{\cite[Theorem 3.26]{KKOP21}}]
For $\varpi \in\wl^+$ and $\mu,\zeta \in \weyl\varpi$ such that $\mu\preceq \zeta$, the determinantial module $M(\mu,\zeta)$ is affreal.
\end{theorem}

\begin{proposition}%
[{\cite{KKKO18,KKOP21}}]\label{prop: l=r} 
Let $M$ and $N$ be simple modules such that one of them is affreal.
Then there exists a unique $R$-module homomorphism $\Rr_{M,N} \in \HOM_R(M,N)$ satisfying
$$\HOM_R(M \conv N,N \conv M)=\bfk\, \Rr_{M,N}.$$
\end{proposition}
We call the homomorphism $\Rr_{M,N}$ the \emph{$R$-matrix}.

\begin{definition} \label{def: inv}
For simple $R$-modules $M$ and $N$ such that one of them is affreal, we define
\begin{align*}
\La(M,N) &\seteq  \deg (\Rr_{M,N}) , \\
\tLa(M,N) &\seteq   \frac{1}{2} \Bigl( \La(M,N) + \bl\wt(M), \wt(N)\br \Bigr) , \\
\de(M,N) &\seteq  \frac{1}{2} \big( \La(M,N) + \La(N,M)\big).
\end{align*}
\end{definition}
It is proved in \cite{KKKO18,KKOP21} that the invariants $\tLa(M,N)$ and $\de(M,N)$ in Definition~\ref{def: inv} belong to $\Z_{\ge 0}$.

\smallskip

For simple modules $M$ and $N$, $M \hconv N$ and  $M \sconv N$ denote the head and the socle of $M \conv N$, respectively.

\Prop[{\cite[Lemma 3.1.4]{KKKO15}, \cite[Corollary 4.1.2]{KKKO18}}]
\label{prop: simple head and socle}
Let $M$ and $N$ be simple $R$-modules such that one of them is affreal. 
\bnum
\item The image of $\rmat{M,N}$ is equal to $M\hconv N$ and
$N\sconv M$.
\item  $M \hconv N$ and  $M \sconv N$ are simple modules
and each of them appears exactly once in the composition series of $M \conv N$
\ro up to a grading shift{\rm)}.
\item Assume that $N$ is affreal.
\bna
\item If a simple subquotient $S$ of $M \conv N$ is not isomorphic to $M\hconv N$, then
$\La(S,N)<\La(M\hconv N,N)=\La(M,N)$.
\item
If  a simple subquotient $S$ of $M \conv N$ is not isomorphic to $M\sconv N$, then
$\La(N,S)<\La(N,M\sconv N)=\La(N,M)$.
\ee
\item  If $M$ and $N$ are self-dual, then $  q^{\tLa(M,N)} M \hconv N$ is a self-dual simple module.
\item \label{it: de 0 simple} The following conditions are equivalent$\colon$
\bna
\item $M \conv N \simeq N \conv M$ up to a grading shift,
\item$M \conv N$ is simple,
\item $\de(M,N)=0$,
\item $M \hconv N \simeq M \sconv N$ up to a grading shift.
\ee
\ee
 \end{proposition}

\begin{proposition}[{\cite[Corollary 3.18]{KKOP21}}] \label{prop: enough big commute}
Let $M$ be an affreal simple module.
Let $X$ be an $R$-module in $R\gmod$. Let $n \in \Z_{>0}$ and assume that any simple
subquotient $S$ of $X$ satisfies $\de(M, S) \le n$. Then any simple subquotient $N$ of
$M \conv X$ satisfies $\de(M, N) < n$. In particular, any simple subquotient of $M^{\circ n} \conv  X$
strongly commutes with $M$.
\end{proposition}

An ordered sequence of simple modules $\uL=(L_1,\ldots,L_r)$ is called \emph{almost affreal} if all $L_i$ $(1 \le i \le r)$ are  affreal  except for at most one.

\begin{definition} An almost affreal sequence $\uL$ of simple modules is called a \emph{normal sequence} if the composition of $R$-matrices
\begin{align*}
\rmat{\uL} & \seteq \prod_{1 \le i < k \le r} \rmat{L_i,L_k }   = (\rmat{L_{r-1},L_r}) \circ \cdots (\rmat{L_2,L_r}\circ \cdots \circ \rmat{L_2,L_3}) \circ (\rmat{L_1,L_r} \circ \cdots \circ \rmat{L_1,L_2}) \\
& \colon \; q^{\sum_{1 \le i < k \le r} \La(L_i,L_k)} L_1 \conv \cdots \conv L_r \To L_r \conv \cdots L_1 \ \  \text{ does not vanish.}
\end{align*}
\end{definition}

\begin{lemma}  [{ \cite[\S 2.3]{KK19}, \cite[\S 2.2]{KKOP22B}}] \label{lem: normality}
Let $\uL$ be an almost affreal sequence of simple modules.
If $\uL$ is normal, then the image of $\rmat{\uL}$ is simple and coincide with the head of $L_1 \conv \cdots \conv L_r$ and also with the socle of $L_r \conv \cdots\conv L_1$, up to grading shifts.
\end{lemma}

\begin{lemma}[{\cite[Proposition 3.2.13.]{KKKO18}}] \label{lem: commute head La}
Let $(A,B,C)$ be an almost affreal sequence. Then we have the followings:
\bnum
\item $\La(A,B\hconv C) =\La(A,B)+\La(A,C)$ if $A$ and $B$ commute.
\item $\La(A\hconv B, C) =\La(A,C)+\La(B,C)$ if $B$ and $C$ commute.
\ee
\end{lemma}

For a given almost affreal sequence $\uL$ of $R$-modules, the sufficient conditions for $\uL$ being normal are studied in \cite{KK19,KKOP22B}. In this paper, we will use the conditions frequently.

\subsection{Commuting families} \label{subsec: commuting family}
Let $J$ be an index set. We say that a family  of affreal simple modules $\calM= \{ M_j \}_{j \in J}$ in $R\gmod$ is a
\emph{commuting family} if
$$M_i \conv M_j \simeq M_j \conv M_i \quad \text{up to a grading shift for any $i,j \in J$.}$$

For a  commuting family
$\calM= \{ M_j \ | \ j \in J\}$ in
$R\gmod$, let us take $\uplambda\cl \Z^{\oplus J} \times \Z^{\oplus J} \to \Z$ such that
\eq
&&
\uplambda(\bfe_i,\bfe_j)-\uplambda(\bfe_j,\bfe_i) = \La(M_i,M_j)
\qt{for any $ i,j \in J$.}\label{def:lambda}
\eneq
Here $\{ \bfe_j \mid j \in J\}$ is the standard basis of $\Z^{\oplus J}$.
Then there exists a family $\{ \calM_\uplambda(\bfa) \mid \bfa \in \Z_{\ge 0}^{\oplus J}\}$ of simple modules in $\scrC_w$ such that
\eq&&\label{def:prod}
\ba{l}
\calM_\uplambda(0)  = \bfone, \qquad
\calM_\uplambda(\bfe_j)  =M_j \quad \text{ for any }j\in J, \\
\calM_\uplambda(\bfa) \conv  \calM_\uplambda(\bfb) \simeq
q^{ -\uplambda(a,b) } \calM_\uplambda(\bfa+\bfb) \quad \text{ for any }\bfa,\bfb \in \Z_{\ge 0}^{J}.
\ea
\eneq
We sometimes omit $_\uplambda$ for notational simplicity.

\Rem Note that $\uplambda(\bfe_i,\bfe_j)=\tLa(M_i,M_j)$ satisfies
condition \eqref{def:lambda}. Moreover,
if all the $M_i$'s are self-dual, then $\calM_\uplambda(\bfa)$
is self-dual for any $\bfa\in \Z_{\ge 0}^{J}$.
We usually take this choice of $\uplambda$.
\enrem
\begin{definition}
A commuting family $\calM=\{ M_j \mid j \in J\}$ is called \emph{\dsc} 
if $\bfa,\bfb \in \Z_{\ge 0}^{\oplus J}$ satisfies
$\calM(\bfa) \simeq q^s \calM(\bfb)$ for some $s \in \Z$, then we have $\bfa=\bfb$.
\end{definition}

The following lemma is obvious.
\begin{lemma}
Let $\calM = \{ M_j \mid j \in J \}$ be a commuting family.
Then  it is \dsc if and only if the set $\big\{ \big[ \calM(\bfa)\big] \mid \bfa \in \Z_{\ge 0}^{\oplus J} \big\}$ in $K(R\gmod)$ is linearly-independent over $\Z[q^{\pm 1}]$.
\end{lemma}

\subsection{Localization of $\scrC_w$} Throughout this subsection, we fix $w \in \weyl$ and set
$$I_w\seteq\stt{i\in I\mid w\varpi_i\not=\varpi_i}.$$ For the notational  simplicity, let us write
$$ C_i \seteq M(w \varpi_i,\varpi_i) \in R\gmod \qquad  \text{ for $i \in I$.}$$
Then $\stt{ \Upomega([C_i]) \mid i \in I_w}$ forms the set of frozen variables of $\calA_\bbA(\n(w))$. For each $\mu = \sum_{i \in I} \mu_i\varpi_i  \in \wl^+$, we set
$C_\mu= M(w\mu,\mu)$, which is a self-dual convolution product $q^c \conv_{i \in I} C_i^{\circ \mu_i}$ for some $c\in \Z$.

It is proved in~\cite{KKOP21,KKOP22B,KKOP23} that there exists a monoidal abelian category $\tscrC_w = \scrC_w[C_i^{\circ -1} \mid  i \in I]$ with a tensor product $\conv$, a degree shift functor $q$
and a monoidal \emph{exact fully faithful} functor $\Phi_w\cl \scrC_w \to \tscrC_w$ satisfying the following properties:
\begin{eqnarray} &&
\parbox{75ex}{
\bnA
\item The objects $\Phi_w(C_i)$ are invertible objects in $\tscrC_w$; that is, there exists an object $\Phi_w(C_i)^{-1}$ in $\tscrC_w$ such that $\Phi(C_i) \conv \Phi_w(C_i)^{-1} \simeq 1$
and $\Phi(C_i)^{-1} \conv \Phi(C_i) \simeq 1$.
\item The category $\tscrC_w$ is universal to $\scrC_w$ in the following sense: For any monoidal functor $\Psi\cl \scrC_w \to \calT$ to another monoidal category $\calT$ in which $\Psi(C_i)$ is invertible for every $i \in I$, there exists a monoidal functor
$\Psi'\cl\tscrC_w \to \calT$ such that $\Psi \simeq \Psi' \circ \Phi_w$. Moreover, $\Psi'$ is unique up to a unique isomorphism.
\item There exists a  commuting family of simple objects $\{ \tC_\mu \ |  \ \mu \in \wl\}$ such that $\tC_\mu \simeq \Phi_w(C_\mu)$ for every $\mu \in \wl^+$ and $\tC_\mu \conv \tC_{\mu'} \simeq q^{\sfH(\mu,\mu')}\tC_{\mu+\mu'}$ for every $\mu,\mu' \in \wl$. Here $\sfH$ denotes the bilinear form on $\wl$ given by $\sfH(\mu,\mu')= (\mu,w\mu'-\mu')$.
\item \label{it: simple} Every simple object in $\tscrC_w$ is isomorphic to $\Phi_w(S) \conv \tC_\mu$ for some simple object $S \in \scrC_w$ and $\mu \in \wl$.
\ee
}\label{eq: tscrC_w}
\end{eqnarray}
(For the precise properties, see loc.\ cit.)

\begin{theorem} [{\cite{KKOP21} (see also \cite[Remark 3.6]{KKOP22B})}]   There exists an $\bbA$-algebra isomorphism
$$\widetilde{\Upomega}  \cl \calK_\bbA(\tscrC_w)  \seteq \bbA \tens_{\Z[q^{\pm1}]} K(\tscrC_w) \isoto \calA_\bbA(\n^w) \quad \text{ such that } \quad \widetilde{\Upomega}|_{\calK_\bbA(\scrC_w)} = \Upomega.$$
Here $K(\tscrC_w)$ denotes the Grothendieck ring of $\tscrC_w$.
\end{theorem}

A pair $(\ve\cl X \otimes Y \to \bfone, \; \eta\cl\bfone \to Y \otimes X)$ of morphisms in a monoidal category with a unit object $\bfone$ is called an \emph{adjunction} if
the following two conditions hold.
\bna
\item The composition $X \simeq X \otimes \bfone \To[X \otimes \eta] X \otimes Y \otimes X \To[\ve \otimes X] \bfone \otimes X \simeq X$ is the identity.
\item  The composition $Y \simeq \bfone \otimes Y \To[\eta \otimes Y] Y \otimes X \otimes Y \To[Y \otimes \ve] Y \otimes \bfone \simeq Y$ is the identity.
\ee
In the case when $(\ve,\eta)$ is an adjunction, we say that $X$ is a \emph{left dual to $Y$}, $Y$ is a \emph{right dual to $X$} and $(X,Y)$ is a \emph{dual pair}.

\begin{theorem} [\cite{KKOP21,KKOP22B}]
The monoidal category $\tscrC_w$ is rigid {\rm ;} i.e.,
every object of $\tscrC_w$ has a right dual and a left dual in $\tscrC_w$.
\end{theorem}

\subsection{Determinantial modules and monoidal clusters} In this subsection, we denote by $\calC$ the category of $\scrC_w$ or $\widetilde{\scrC}_w$.
Recall that
$$\parbox{75ex}{$\scrA=\calK_\bbA(\calC)$  has a quantum cluster algebra structure via an isomorphism $\Omega =\Upomega$ or $\widetilde{\Upomega}$.}$$

Let $\ii=(i_1,\ldots,i_r)$ be a reduced sequence of $w \in \weyl$.
For $k$ such that $1\le k\le r$, set
$$M_k^\ii=M(w^\ii_{\le k}\varpi_{i_k},\varpi_{i_k})$$
(see Proposition~\ref{prop: determinnatial} for the notation).

\begin{proposition} [{\cite[Theorem 4.12]{KKOP18}}] \label{prop: La =l } Let $\ii=(i_1,\ldots,i_r)$ be a reduced sequence of $w \in \weyl$.
For $ s < t \in \K$, we have
$$  -\La(M_s^\ii,M_t^\ii) = (\varpi_{i_s}-w^\ii_{\le s}\varpi_{i_s}, \varpi_{i_t}+w^\ii_{\le t}\varpi_{i_t}  )  = l^\ii_{st} = (L^\ii)_{st}.$$
\end{proposition}

We say that a commuting family $\calM=\stt{M_i}_{i\in \K}$ in $\calC$ is a \emph{monoidal cluster} if there exists
a quantum seed  $(\{ X_i \}_{i \in \K}, L=(l_{i,j})_{i,j \in \K}, \tB=(b_{i,j})_{i \in \K,j\in \Kex} )$ of $\scrA$ such that
$$   X_i = \Omega( q^{\frac{1}{4}(\wt(M_i),\wt(M_i))}[M_i]) \quad \text{ and } \quad l_{i,j} = -\La(M_i,M_j).$$

Note that every monoidal seed is  \dsc \ since the quantum cluster monomials in a cluster are linearly independent  over $\bbA$ by the definition.

With Proposition~\ref{prop: determinnatial} and \eqref{eq: seed ii}, Proposition~\ref{prop: La =l }  tells that
\begin{align} \label{eq: calM ii}
    \calM^\ii \seteq \{ M^\ii_k  \}_{1 \le k \le r} \text{ is a monoidal cluster in $\scrC_w$,}
\end{align}
for any reduced sequence $\ii=(i_1,\ldots,i_r)$ of $w$.
We call $\calM^\ii$ the \emph{GLS-cluster} (associated with $\ii$).

\section{Quasi-Laurent family and Laurent family} \label{sec: Laurent family}

In this section, we introduce the notion of \ql families and
Laurent families, which allows us to associate two vectors in $\Z^{J}$ with each simple module.

\subsection{Definition} Let $J$ be a finite index set.
Let $\scrC$ be a full monoidal subcategory of $R\gmod$ stable by taking
subquotients, extensions and grading shifts.

\begin{lemma}\label{lem:well-defined}
Let $\calM=\{ M_j \mid j \in J\}$ be an {\em \dsc} commuting family of affreal simple objects in $\scrC$
and $X$ a simple module in $\scrC$.
\bnum
\item \label{it: i} If $\bfa$, $\bfa'$, $\bfb$, $\bfb'\in\Z_{\ge 0}^{J}$
satisfy
$ X \hconv  \calM(\bfa)  \simeq \calM(\bfb)$ and $X \hconv \calM(\bfa')    \simeq \calM(\bfb')$
up to grading shifts, then one has
$\bfb-\bfa=\bfb'-\bfa'$.
\item\label{it: ii} If
$\bfa$, $\bfa'$, $\bfb$, $\bfb'\in\Z_{\ge 0}^{J}$
satisfy
$ \calM(\bfa) \hconv X \simeq \calM(\bfb)$ and $ \calM(\bfa')  \hconv X  \simeq \calM(\bfb')$
up to grading shifts, then one has
$\bfb-\bfa=\bfb'-\bfa'$.
\ee
\end{lemma}
\begin{proof}
Since the proof are similar, we prove only~\eqref{it: i}.
We have
\eqn
\calM(\bfb+\bfa') \simeq \calM(\bfb) \hconv \calM(\bfa')
&&\simeq( X \hconv \calM(\bfa) )  \hconv  \calM(\bfa')
\simeq \head(X \conv  \calM(\bfa') \conv  \calM(\bfa)  )\\
&&\simeq( X \hconv \calM(\bfa') )  \hconv  \calM(\bfa)
\simeq \calM(\bfb'+\bfa),
\eneqn
and hence we have $\bfa'+\bfb=\bfa+\bfb'$ since $\catC$ is \dsc.
\end{proof}

\begin{definition} \label{def: quasi-Laurent}
We say that a commuting family $\calM= \{ M_j \mid j \in J \}$
of affreal simple objects of $\scrC$ is a \emph{\ql family} in $\scrC$  if
$\calM$ satisfies the following conditions:
\bna
\item\label{it:dsc}
$\calM$ is \dsc, and
\item if a simple module $X$ commutes with all $M_j$ $(j \in J)$,
then there exist $\bfa,\bfb \in \Z_{\ge 0}^{J}$ such that
$$
X \conv \calM(\bfa) \simeq \calM(\bfb)\qt{up to a grading shift}.
$$
\setcounter{myc}{\value{enumi}}
\ee
If $\calM$ satisfies \eqref{it:dsc} and \eqref{it:Laur} below, then we say that $\calM$
is a {\em Laurent family}:
\bna
\setcounter{enumi}{\value{myc}}
\item  if a simple module $X$ commutes with all $M_j$ $(j \in J)$,
then there exists $\bfb\in \Z_{\ge 0}^{J}$ such that $X\simeq\calM(\bfb)$.
\label{it:Laur}
\ee
Note that a Laurent family is a quasi-Laurent family.
\end{definition}

\begin{lemma} \label{lem: a,a',b,b'}
Let $\calM = \{ M_j \mid j \in J\}$ be a \ql family in $\scrC$. Then we have the followings$\colon$
\bnum
\item for any simple module $X \in \scrC$, there exist $\bfa,\bfb \in \Z^J_{\ge 0}$ such that $X \hconv \calM(\bfa) \simeq \calM(\bfb)$ up to a grading shift,
\item for any simple module $X \in \scrC$, there exist $\bfa',\bfb' \in \Z^J_{\ge 0}$ such that $X \sconv \calM(\bfa') \simeq \calM(\bfb')$
up to a grading shift.
\ee
\end{lemma}

\begin{proof}
Since the proofs are similar, we shall only prove the first statement. Let take $\bfa^{(1)} \in \Z^J_{\ge 0}$ such that
$\bfa^{(1)}_j \gg 0$ for all $j \in J$. Then Proposition~\ref{prop: enough big commute} tells that the simple module $Y \seteq X \hconv \calM(\bfa^{(1)})$
commutes with all $M_j$. Since $\calM$ is a \ql family, there exists $\bfa^{(2)},\bfb \in \Z_{\ge0}^J$ such that $Y \conv \calM(\bfa^{(2)}) \simeq \calM(\bfb)$.
Hence, by taking $\bfa = \bfa^{(1)}+\bfa^{(2)}$, we have
$$ X  \hconv \calM(\bfa) \simeq \calM(\bfb),$$
as desired.
\end{proof}

By applying the similar argument in the lemma above to composition factors of $X \conv \calM(\bfa)$, we have the following corollary:

\begin{corollary}
Let $\calM $ be a \ql family in $\scrC$.
Then, for any $X\in \scrC$, there exist $\bfa \in \Z_{\ge 0}$, a finite index set $\sfS$, $c(s) \in \Z$ and $\bfb(s) \in \Z_{\ge 0}^J$ $(s \in \sfS)$ such that
\begin{align} \label{eq: Laurent}
[X \conv \calM(\bfa)] = \sum_{s \in \sfS} q^{c(s)}[\calM(\bfb(s))].
\end{align}
\end{corollary}

\begin{remark} \label{rem: Laurent phenomenon}
The above corollary tells that for every module $X$ in $\scrC$ and a \ql family  $\calM = \{ M_j \mid j \in J\}$,
the isomorphism class $[X]$ of $X$ in $K(\scrC)$
can be expressed as an element in the Laurent polynomial ring
$\Z[q^{\pm 1}][[M_j]^{\pm1} \mid j \in J]$ with positive coefficients.
\end{remark}

\begin{proposition} \label{prop:Layrent}
Assume that $\calK_\bbA(\catC)$ has a quantum cluster algebra structure, and let
$\calM=\{ M_k \mid  k\in \K\}$ be a monoidal cluster in $\scrC$. Then the commuting family $\calM$
is a quasi-Laurent family.
In particular,
the isomorphism class $[X]$ of $X$ in $K(\scrC)$
can be expressed as an element in the Laurent polynomial $\Z[q^{\pm 1}][\;[M_j] \mid j \in J]$ with positive coefficients.

If moreover every $[M_k]$ is prime in $K(\catC)\vert_{q=1}$,
then $\calM$ is a Laurent family.
\end{proposition}

Note that if $K(\catC)\vert_{q=1}$ is factorial, then every $[M_k]$ is prime in $K(\catC)\vert_{q=1}$ (see \cite{GLS13A}).

\begin{proof}
Let $X$ be a simple module in $\scrC$ commuting with all $M_k$ $(k\in \K)$. The quantum Laurent phenomenon states that
there exists $\bfa \in \Z_{\ge 0}^\K$ such that
\begin{align} \label{eq: graded decomposition}
[X \conv \calM(\bfa)] = [X]\cdot [\calM(\bfa)] = \sum_{s \in \sfS} c(s)[\calM(\bfb^{(s)})]
\end{align}
for some $\ell \in \Z_{\ge 1}$, $\bfb^{(s)} \in \Z_{\ge 0}^\K$  and $c(s) \in \Z[q^{\pm 1/2}]$.  Since
$X \conv \calM(\bfa)$ is simple, the right-hand side of~\eqref{eq: graded decomposition} must coincide with $q^c\calM(\bfc)$ for some $\bfc \in \Z_{\ge 0}^J$ and $c \in \Z$.
Hence $\catC$ is a \ql family.

\smallskip
Let us show that $\calM$ is a Laurent family.
If $X\conv \calM(\bfa)\simeq\calM(\bfb)$, then
we have $\bfa_k\le\bfb_k$ for all $k$, since
each $[M_k]$ is a prime element of $K(\catC)\vert_{q=1}$.
Hence setting $\bfc=\bfb-\bfa\in\Z_{\ge0}^\K$,
we have
$X\conv\calM(\bfa)\simeq\calM(\bfc)\conv \calM(\bfa)$,
which implies that
$X\simeq\calM(\bfc)$.
Hence $\calM$ is a Laurent family.
\end{proof}

\begin{definition} \label{def: g-vectors}
For a simple module $X \in \scrC$ and a \ql family $\calM=\{ M_j \mid j \in J \}$, take $\bfa,\bfb,\bfa',\bfb' \in \Z_{\ge 0}^J$ such that
$X \hconv \calM(\bfa) \simeq \calM(\bfb)$ and $X \sconv \calM(\bfa') \simeq \calM(\bfb')$ up to grading shifts. Then we define
$$
\bfg_\calM^R(X) \seteq \bfb -\bfa \quad \text{and} \quad \bfg_\calM^L(X) \seteq \bfb' -\bfa' \in \Z^J.
$$
\end{definition}

\begin{remark} \hfill
\ben
\item For a \ql family $\calM$ in $\scrC$, $\bfg_\calM^R$ and $\bfg_\calM^L$ are well-defined by Lemma~\ref{lem:well-defined}.
\item For a reduced sequence $\ii$ of $w$ and its \ql family $\calM^\ii$, we write $\bfg_\ii^R$ and $\bfg_\ii^L$ instead of $\bfg_{\calM^\ii}^R$
and $\bfg_{\calM^\ii}^L$, respectively.
\item The map $\bfg_\ii^R$ and $\bfg_\ii^L$ for the \ql family $\calM^\ii$ can be extended to the set $\Irr(\tscrC_w)$ of the isomorphism classes of self-dual simples in $\tscrC_w$.
\ee
\end{remark}

The following lemma can be proved by the same arguments in \cite{KK19}.

\begin{lemma}\label{lem:gvec}
 Let  $\calM$ be a \ql family in $\scrC$
and $X$ a simple module in $\scrC$.
\bnum
\item If $\bfa,\bfb \in \Z_{\ge 0}^{J}$ satisfy $\bfb-\bfa = \bfg_\calM^R(X)$ $($resp.\  $\bfb-\bfa = \bfg_\calM^L(X))$, then we have
$$ \text{$X \hconv \calM(\bfa) \simeq \calM(\bfb)$
\ro resp.\ $X\sconv \calM(\bfa) \simeq \calM(\bfb)${\rm)} up to a grading shift.}$$
\item For any $\bfa \in \Z_{\ge 0}^{J}$, we have
$$  \bfg_\calM^R(X \hconv \calM(\bfa)) =  \bfg_\calM^R(X)+\bfa \quad \text{and}\quad  \bfg_\calM^L(X \sconv \calM(\bfa)) =  \bfg_\calM^L(X)+\bfa.$$
\item The maps $ \bfg_\calM^R$ and $ \bfg_\calM^L$ from $\Irr(\scrC)$ to $\Z^{J}$ are injective.
\ee
\end{lemma}

\section{PBW-decomposition vector and GLS-seed} \label{sec: PBW and GLS}
In this section, we recall the PBW-basis,
and we state the relations of the $g$-vector and the PBW-decomposition vector.

\subsection{PBW-decomposition vector} Let us take
$w \in \weyl$ and its reduced sequence $\ii=(i_1,\ldots,i_r)$.
Recall that
$$\K=[1,r]=\Kex  \sqcup\, \Kfr, \ \Kex=\stt{k\in\K\mid k_+\le r},
\  \be^\ii_k \in \Delta^+ \cap w\Delta^-,  \ S^\ii_k  \text{ and }   M^\ii_k \   (1\le k \le r),$$
and note that
\begin{align} \label{eq: M and S}
 M^\ii_k \simeq S^\ii_k \hconv M^\ii_{k_-}.
\end{align}
Then $\calM^\ii \seteq \stt{M_k^\ii\mid k\in \K}$ is a commuting family.
For any $\bfa = (\bfa_k)_{1 \le k \le r} \in \Z_{\ge 0}^\K$, the convolution product
$$ P_\ii(\bfa) \seteq q^{\frac{1}{2}  \sum_{k=1}^r \bfa_k(\bfa_k-1) \sfd_{i_k}  } {S^\ii_r}^{\;\circ \bfa_r} \conv \cdots \conv {S^\ii_1}^{\;\circ \bfa_1}  $$
has a self-dual simple head. Conversely, every self-dual simple module in $\scrC_w$ is isomorphic to $\head(P_\ii(\bfa))$ for some $\bfa \in \Z_{\ge0}^\K$ in a unique way ( \cite{McNa15,TW16}).
We call $\{ P_\ii(\bfa) \mid \bfa \in \Z_{\ge0}^\K \}$
the \emph{PBW-basis of $\scrC_w$ associated with $\ii$}.

For a simple module $X$ such that $X\simeq\head(P_\ii(\bfa))$, we set
$$\PBW_\ii(X) \seteq \bfa = (\bfa_1,\ldots,\bfa_r).$$

The following lemma tells that the operation
$\PBW_\ii( - \hconv \calM^\ii(\bfa))$ on the set of simple modules behaves very nicely, where $\calM^\ii(\bfa)$ is defined in \eqref{def:prod}.

\begin{lemma}   [{cf.\ \cite[Lemma 3.11 and Proposition 3.14]{KK19}}]
\label{lem: PBW addition}
For $M=\calM^\ii(\bfa)$ with $\bfa \in \Z_{\ge0}^{\K}$ and a simple module $X$, we have
$$
\PBW_\ii(X \hconv M) =\PBW_\ii(X)+\PBW_\ii(M).
$$
In particular $\bfc\seteq\PBW_\ii\bl\calM^\ii(\bfa)\br$ is given by
$\bfc_k=\sum_{j}\bfa_j$
where $j$ ranges over $j\in[1,r]$ such that $j\ge k$ and
$\ii_j=\ii_k$.
\end{lemma}
\begin{proof}
It is enough to show it when $M=M^\ii_k$. Note that
$$  X \simeq \head\left( \conv_{1 \le k \le r}^{\to} {S^\ii_k}^{\circ n_k} \right) = {S^\ii_{r}}^{\circ n_r } \hconv Y$$
where $\bfn=(\bfn_1,\ldots,\bfn_r) = \PBW_\ii(X)$,
$Y \simeq \head\left( \displaystyle\conv_{1 \le k \le r-1}^{\to} {S^\ii_k}^{\circ n_k} \right)$ \and $\displaystyle\conv_{p \le k \le q}^{\to} X_k$  denotes the ordered convolution product $X_q \conv X_{q-1} \conv \cdots X_{p+1} \conv X_p$ for $X_k$'s in $R\gmod$.

\noindent
If  $r > k$, then $({S^\ii_r}^{\circ n_r},M^\ii_k)$ is unmixed and hence
${S^\ii_r}^{\circ n_r} \conv Y \conv M^\ii_k$ has a simple head. We have
\begin{align*}
X \hconv M^\ii_k &  \simeq \head({S^\ii_r}^{\circ n_r} \conv Y \conv M^\ii_k) \allowdisplaybreaks \\
& \simeq {S^\ii_r}^{\circ n_r} \hconv (Y \hconv M^\ii_k) \simeq {S^\ii_r}^{\circ n_r} \hconv \head\left( \displaystyle\conv_{1 \le k \le r-1}^{\to} {S^\ii_k}^{\circ \bfc_k} \right),
\end{align*}
where $\bfc =\PBW_\ii(Y)+\PBW_\ii(M^\ii_k)$ by induction on $r$. Thus our assertion follows in this case.

\noindent
If  $r = k$,  then $M^\ii_r $ commutes with all the objects of $\catC_w$,
and hence we have
\begin{align*}
X \hconv M^\ii_r &  \simeq \head({S^\ii_r}^{\circ n_r} \conv M^\ii_r \conv {S^\ii_{r-1}}^{\circ n_{r-1}} \conv \cdots {S^\ii_{1}}^{\circ n_{1}} ) \allowdisplaybreaks\\
& \simeq \head({S^\ii_r}^{\circ n_r} \conv M^\ii_r) \hconv Y   \simeq \head({S^\ii_r}^{\circ n_r} \conv (S^\ii_r \hconv M^\ii_{r_-})) \hconv Y \allowdisplaybreaks \\
& \simeq \head({S^\ii_r}^{\circ n_r+1}  \hconv M^\ii_{r_-} ) \hconv Y   \simeq {S^\ii_r}^{\circ n_r+1}  \hconv (M^\ii_{r_-}  \conv Y ),
 \end{align*}
where the last isomorphism follows from the commutativity of
$M^\ii_{r_-}$ and $Y$.
Then our assertion follows from the induction hypothesis.
\end{proof}

The lemma above gives a direct proof of the following corollary
although it follows immediately from Proposition~\ref{prop:Layrent}.
\Cor
For any reduced sequence $\ii$ for $w$, the commuting family $\calM^\ii$ is
\dsc.
\encor

The following proposition is proved in \cite[Proposition 2.11]{KK19} for symmetric quiver Hecke algebra and the same proof also works for the general case:

\begin{proposition} \label{prop: normal pbw}
For any $\bfa = (\bfa_1,\ldots,\bfa_r)$, $\bfb = (\bfb_1,\ldots,\bfb_r) \in \Z_{\ge 0}^\K$, the ordered sequence
$$
({S^\ii_r}^{\circ \bfa_r},(S_{r-1}^\ii)^{\circ \bfa_{r-1}},\ldots, {S^\ii_1}^{\circ \bfa_1},{M^\ii_1}^{\circ\bfb_1},{M^\ii_2}^{\circ\bfb_2},\ldots,{M^\ii_r}^{\circ\bfb_r})
$$
is a normal sequence.
\end{proposition}

The statement and proof of following proposition are the same as \cite[Proposition 3.14]{KK19} even though
\cite{KK19} dealt only with symmetric quiver Hecke algebras. Here we repeat
it in order to show relations between explicit $\Z_{\ge 0}$-vectors associated
with a simple module $X$ in $\scrC_w$ for the readers' convenience.

\begin{proposition} \label{prop: existence a,b}
For a simple module $X$ in $\scrC_w$ or $\tscrC_w$, there exist  $\bfa, \bfb \in \Z_{\ge 0}^\K$ such that
$$
X \hconv \calM^\ii(\bfa) \simeq \calM^\ii(\bfb)  \qquad \text{ up to a grading shift.}
$$
\end{proposition}

\begin{proof} In this proof,  we sometimes drop $^\ii$ for notational simplicity.
It is enough to consider when $X \in \scrC_w$ by~\eqref{it: simple} in~\eqref{eq: tscrC_w}.
Note that there exists a unique $\bfc = (\bfc_1,\ldots,\bfc_r) \in \Z^\K_{\ge 0}$ such that
$$   X \simeq \head(P_\ii(\bfc)) \simeq \head({S^\ii_r}^{\circ \bfc_r} \conv \cdots \conv  {S^\ii_1}^{\circ \bfc_1}) .$$

Set $\bfc_+ \seteq \sum_{k=1}^r  \bfc_{k_+} \bfe_k= \sum_{j \in \K} \bfc_k \bfe_{k_-}$, where $\{ \bfe_j \mid j \in \K\}$ is the standard basis of $\Z^\K$ such that $\bfc = \sum_{j \in \K} \bfc_j \bfe_j$.   Then
we have $\calM^\ii(\bfc_+) \simeq M_{1_-}^{\circ \bfc_1} \conv \cdots \conv M_{r_-}^{\circ \bfc_r}$. By~\eqref{eq: M and S} and Proposition~\ref{prop: normal pbw}, we have
\begin{align*}
& X \hconv \calM^\ii(\bfc_+) \simeq \head(S_r^{\circ \bfc_r}  \conv  S_{r-1}^{\circ \bfc_{r-1}}  \conv \cdots \conv S_1^{\circ \bfc_{1}} \conv \calM^\ii(\bfc_+)  )\allowdisplaybreaks \\
&\simeq \head(S_r^{\circ \bfc_r}  \conv  S_{r-1}^{\circ \bfc_{r-1}}  \conv \cdots \conv S_2^{\circ \bfc_{2}} \conv S_1^{\circ \bfc_{1}} \conv M_{1_-}^{\circ \bfc_1} \conv M_{2_-}^{\circ \bfc_2} \conv \cdots \conv M_{r_-}^{\circ \bfc_r} )\allowdisplaybreaks \\
&\simeq \head((S_r^{\circ \bfc_r}  \conv  S_{r-1}^{\circ \bfc_{r-1}}  \conv \cdots \conv S_2^{\circ \bfc_{2}}) \conv (S_1^{\circ \bfc_{1}} \hconv M_{1_-}^{\circ \bfc_1}) \conv (M_{2_-}^{\circ \bfc_2} \conv \cdots \conv M_{r_-}^{\circ \bfc_r} ))\allowdisplaybreaks \\
&\simeq \head((S_r^{\circ \bfc_r}  \conv  S_{r-1}^{\circ \bfc_{r-1}}  \conv \cdots \conv S_1^{\circ \bfc_{2}}) \conv  M_{1}^{\circ \bfc_1} \conv (M_{2_-}^{\circ \bfc_2} \conv \cdots \conv M_{r_-}^{\circ \bfc_r} )) \allowdisplaybreaks\\
&\simeq \head((S_r^{\circ \bfc_r}  \conv  S_{r-1}^{\circ \bfc_{r-1}}  \conv \cdots \conv S_2^{\circ \bfc_{2}})  \conv (M_{2_-}^{\circ \bfc_2} \conv \cdots \conv M_{r_-}^{\circ \bfc_r} )\conv  M_{1}^{\circ \bfc_1} )\allowdisplaybreaks \\
& \simeq \cdots  \simeq \head(M_r^{\circ \bfc_r} \conv \cdots \conv M_1^{\circ \bfc_1}) \simeq \calM^\ii(\bfc),
\end{align*}
which implies our assertion.
\end{proof}
As seen by the proof of the above proposition and Proposition~\ref{prop:Layrent}, we have
\begin{proposition} \label{prop: Laurent commuting}
The commuting family $\calM^\ii$ is a Laurent family.
Moreover for a simple module
$M$, two  vectors  $\bfa=\PBW_\ii(M)$ and $\bfg\seteq\bfg_\ii^R(M)$
are related by
$$\bfg_k=\bfa_k-\bfa_{k_+},\quad \bfa_k=\sum_{j; j\ge k,\ \ii_j=\ii_k}\bfg_j,$$
where $\bfa_{r+1}=0$.
\end{proposition}

The following corollary  can be proved by the same arguments in \cite{KK19}.

\begin{corollary} \label{lem: bijective}
Let $\ii$ be a reduced sequence of $w$.
\ben
\item For a dual pair of simples $(L,R)$ in  $\tscrC_w$, we have
$$ \bfg_\ii^R(L)+ \bfg_\ii^L (R)=0$$
\item The maps $\bfg_\ii^R, \bfg_\ii^L \colon \Irr(\tscrC_w)  \to \Z^\K$  are bijective.
\ee
\end{corollary}

\section{Skew symmetric pairings} \label{sec: g-vector}

In this section, we study skew symmetric pairings induced by the $\Z$-vectors associated with simple modules.

\subsection{Skew-symmetric pairing associated with a \ql family}
Let $\scrC$ be a full monoidal subcategory of $R\gmod$ stable by taking
subquotients, extensions and grading shifts,
and let $\calM=\{ M_j \mid j \in J\}$  be a quasi-Laurent family in $\scrC$ labeled by a finite index set $J$.

For
$X,Y \in \Irr(\scrC)$, let us define
\eq
&&\ba{ll}
&\rmG^R_\calM(X,Y) \seteq \sum_{  a,b \in J} (\bfg^R_\calM(X))_a (\bfg^R_\calM(Y))_b \, \La(M_a,M_b) \qt{and}\\[2ex]
&\rmG^L_\calM(X,Y) \seteq \sum_{  a,b \in J} (\bfg^L_\calM(X))_a (\bfg^L_\calM(Y))_b \, \La(M_a,M_b).
\ea
\eneq

The following lemma immediately follows from Lemma~\ref{lem:gvec}.
\Lemma\label{lem:Grconv}
 For $M=\calM(\bfa)$ with $\bfa \in \Z_{\ge 0}^{J}$ and $X,Y \in \Irr(\scrC)$, we have
\begin{align*}
\rmG^R_\calM(X \hconv M,Y) = \rmG^R_\calM(X,Y)+ \rmG^R_\calM(M,Y),  && \rmG^R_\calM(X,Y) = -\rmG^R_\calM(Y,X), \\
\rmG^L_\calM(X \sconv M,Y) = \rmG^L_\calM(X,Y)+ \rmG^L_\calM(M,Y),  && \rmG^L_\calM(X,Y) = -\rmG^L_\calM(Y,X).
\end{align*}
\enlemma

\begin{proposition} \label{prop: G Y MCC}
Let $X$ be a simple module in $\scrC$.
Then for any $\bfc\in\Z_{\ge0}^J$, we have
$$ \rm{(i)} \;  \La(X,\calM(\bfc)) =\rmG^R_\calM(X,\calM(\bfc)) \quad \text{ and } \qquad \rm{(ii)} \; \La(\calM(\bfc),X) =\rmG^L_\calM(\calM(\bfc),X) .$$
\end{proposition}

\begin{proof} If $X$ is also of the form $\calM(\bfd)$, it is obvious.
Set $Y=\calM(\bfc)$.

\smallskip \noindent
(i)  Note that there exist $\bfa,\bfb \in \Z_{\ge 0}^{J}$ such that $X \hconv \calM(\bfa) \simeq \calM(\bfb)$. Then we have
\begin{align*}
 \La(X,Y)+\La(\calM(\bfa),Y) = \La(X\hconv \calM(\bfa),Y) & = \rmG^R_\calM(X\hconv \calM(\bfa),Y) \\
& = \rmG^R_\calM(X,Y) + \rmG^R_\calM(\calM(\bfa),Y).
\end{align*}
Since $\La(\calM(\bfa),Y)  = \rmG^R_\calM(\calM(\bfa),Y) $, our assertion follows.

\smallskip \noindent
(ii) Similarly, there exist $\bfa,\bfb \in \Z_{\ge 0}^{J}$ such that $\calM(\bfa) \hconv X \simeq \calM(\bfb)$. Then we have
\begin{align*}
\La(Y,\calM(\bfa)) + \La(Y,X) = \La(Y,\calM(\bfa)\hconv X) & = \rmG^L_\calM(Y,\calM(\bfa)\hconv X) \\
& = \rmG^L_\calM(Y,\calM(\bfa)) + \rmG^L_\calM(Y,X).
\end{align*}
Then our assertion follows from the fact that $\La(Y,\calM(\bfa))  = \rmG^L_\calM(Y,\calM(\bfa))$.
\end{proof}

\begin{proposition} \label{prop: lower bound}
For any simple modules $X,Y$ in $\scrC$ such that one of them is affreal, we have
$$    \rmG^R_\calM(X,Y),  \rmG^L_\calM(X,Y) \le \La(X,Y). $$
\end{proposition}

\begin{proof} Since the proofs are similar, we will consider the case of $\rmG^R_\calM$.
Take $\bfa,\bfb\in\Z_{\ge0}^J$ such that $Y \hconv \calM(\bfa) \simeq \calM(\bfb)$. Then we have
\begin{align*}
\rmG^R_\calM(X,Y) + \rmG^R_\calM(X,\calM(\bfa))& = \rmG^R_\calM(X,Y\hconv\calM(\bfa))  \\
& = \La(X,Y \hconv \calM(\bfa)) \le\La(X,Y)+ \La(X,\calM(\bfa))  \\
& = \La(X,Y)+ \rmG^R_\calM(X,\calM(\bfa)) ,
\end{align*}
which yields our assertion.
\end{proof}

\begin{proposition} \label{prop: Gms Lambda}
If simple modules $X$ and $Y$ in $\scrC$ commute and one of them is affreal, then we have
$$   \La(X,Y) = \rmG^R_\calM(X,Y) = \rmG^L_\calM(X,Y).$$
\end{proposition}

\begin{proof} Since the proofs are similar, we will only give the proof for $\rmG^R_\calM$.
By the preceding lemma, we have
\eqn
&&0=\bl\La(X,Y)+\La(Y,X)\br-\bl\rmG^R_\calM(X,Y)+\rmG^R_\calM(Y,X)\br\\
&&\hs{5ex}=\bl\La(X,Y)-\rmG^R_\calM(X,Y)\br+\bl \La(Y,X)-\rmG^R_\calM(Y,X)\br
\ge0,
\eneqn
which implies $\La(X,Y)-\rmG^R_\calM(X,Y)=0$.
\end{proof}

\begin{remark} \label{rmk: bfg depend}
The two invariants $\GR(X,Y)$ and $\GL(X,Y)$ are different in general and depend on the choice of $\calM$.

Let $w_0$ be the longest element of finite type $A_2$.
For a reduced sequence $\ii=(1,2,1)$ of $w_0$,
we have
\begin{align} \label{eq: Lii}
\{ S_1^{\ii}=\lan 1\ran ,  S_2^{\ii}=\lan 12\ran  , S_3^{\ii}=\lan 2\ran \}  \text{ and } \calM^\ii =\{ M_1^{\ii}=\lan 1\ran ,  M_2^{\ii}=\lan 12\ran  , M_3^{\ii}=\lan 21\ran \},
\end{align}
while
\begin{align} \label{eq: Ljj}
\{ S_1^{\jj}=\lan 2\ran ,  S_2^{\jj}=\lan 21\ran  , S_3^{\jj}=\lan 1\ran \}  \text{ and } \calM^\jj =\{ M_1^{\jj}=\lan 2\ran ,  M_2^{\jj}=\lan 21\ran  , M_3^{\jj}=\lan 12\ran \}
\end{align}
for the other reduced sequence $\jj=(2,1,2)$ of $w_0$.
Here $\lan k\ran$ $(k=1,2)$ is a $1$-dimensional $R(\al_k)$-module, and
$\lan 12 \ran$ and  $\lan 21 \ran$ are $1$-dimensional $R(\al_1+\al_2)$-modules (see~\cite{KKK1} for more details on these modules).

Since $A_2$ is symmetric, $\calM' \seteq \mu_1(\calM^\ii)$ is also a Laurent family given as follows:
$$ \calM' = \{ M'_1 =\mu_1(M^\ii_1) \simeq \lan 2\ran, M'_2 = M^\ii_2 \simeq \lan 12 \ran,  M'_3 = M^\ii_3 \simeq \lan 21 \ran\}.$$
Note that $\calM'=\calM^\jj$ (up to an index permutation).
We have
\begin{align*}
& \bfg^R_{\calM^\ii}(\lan 1 \ran) =  \bfg^L_{\calM^\ii}(\lan 1 \ran) = (1,0,0),    \quad \bfg^R_{\calM^\ii}(\lan 2 \ran ) =  (-1,0,1),  \quad
\bfg^L_{\calM^\ii}(\lan 2 \ran ) =  (-1,1,0), \ \   \\
& \bfg^R_{\calM'}(\lan 1 \ran)  = (-1,1,0), \quad    \bfg^L_{\calM'}(\lan 1 \ran)  = (-1,0,1), \quad \bfg^R_{\calM'}(\lan 2 \ran) =
\bfg^L_{\calM'}(\lan 2 \ran) = (1,0,0),
\end{align*}
and
\begin{align*}
&  \La(\lan1\ran,\lan2\ran)  =   \La(\lan2\ran,\lan1\ran)=1,  \qquad \La(\lan1\ran,\lan12\ran)=1, \\
&   \La(\lan1\ran,\lan21\ran) =-1, \quad  \La(\lan12\ran,\lan2\ran) = 1, \quad  \La(\lan21\ran,\lan2\ran) = -1.
\end{align*}
Note that $\La(X,X)=0$ for an affreal simple module $X$. Thus we have
\begin{align*}
& \rmG^R_{\calM^\ii}(\lan1\ran,\lan2\ran) = 1 \times \La(\lan1\ran,\lan21\ran) =-1,  &&\rmG^R_{\calM'}(\lan1\ran,\lan2\ran) =1 \times \La(\lan 12 \ran,\lan 2 \ran) =1, \\
& \rmG^L_{\calM^\ii}(\lan1\ran,\lan2\ran) = 1 \times \La(\lan1\ran,\lan12\ran) =1,
&&\rmG^L_{\calM'}(\lan1\ran,\lan2\ran) =1 \times \La(\lan 21 \ran,\lan 2 \ran) =-1.
\end{align*}
Thus, for a non-commuting pair of simple modules $(X,Y)$ in $\scrC$,
the $\Z$-values $\rmG^R_\calM(X,Y)$ and $\rmG^L_\calM(X,Y)$
\emph{do depend on} the choice of a \ql commuting family $\calM$.

\end{remark}

\subsection{Skew symmetric pairing associated with the GLS-cluster}
Let $w \in \weyl$ and $\ii=(i_1,\ldots,i_r)$ a reduced sequence of $w$.
Let $\calM^\ii$ be the associated GLS-cluster.
For such a Laurent family, we can define $\rmG_{\calM^\ii}^R$ in  terms of PBW-decompositions.

We define a skew-symmetric $\Z$-valued map
$\la^\ii\cl[1,r]\times[1,r]\to\Z$ by
\begin{align} \label{eq: skew form}
\la^\ii _{a,b} \seteq (-1)^{\delta(a>b)} \delta(a \ne b) (\be^\ii_a,\be^\ii_b)
\end{align}
for $1 \le a,b\le r$.

\begin{remark}
The skew symmetric map $\la^\ii$ in~\eqref{eq: skew form} is known when $\g$ is of finite type and $\ii$ is \emph{adapted to a $Q$-datum} (\cite[Proposition 3.2]{HL15},   \cite[Proposition 5.21]{FO21} and \cite[Theorem 5.4]{KO22}).
\end{remark}

Let us recall the notion of $\ii$-box and an affreal simple module $M^\ii[a,b]$ in $\scrC_w$ for an $\ii$-box $[a,b]$, which are introduced in~\cite{KKOP2}.

\bna
\item  For $1 \le a \le  b \le r$ such that $i_a=i_b$, we call an interval  $[a,b]$ an \emph{$\ii$-box}.
\item For an $\ii$-box $[a,b]$, we set $[a,b]_\ii \seteq  \{ u \mid  a \le u \le b, \; i_a=i_u\}$.
\item  For an $\ii$-box $[a,b]$, we set
\begin{align*} M^\ii[a,b] \seteq
M(w^\ii_{\le b}\varpi_{i_a},w^\ii_{<a}\varpi_{i_a})& \simeq  \head\left( \conv_{u \in [a,b]_\ii}^{\to} S^\ii_u \right) \\
&\simeq  S^\ii_b \hconv  M^\ii[a,b_-]   \simeq  M^\ii[a_+,b] \hconv S^\ii_a,
\end{align*}
up to grading shifts. In particular, $M^\ii_k = M^\ii[k_{\min},k]$ and $S^\ii_k = M^\ii[k,k]$.
\ee
Note that $M^\ii[a,b]$ is an affreal simple module in $\scrC_w$.


\begin{proposition} \label{prop: Mi Li}
For $\ii$-boxes $[x,y]$ and $[x',y']$ in an interval $[1,r]$, assume that
\begin{align} \label{eq: left pre-commutative}
 {\rm (a)} \ \  x > x'_- \qquad \text{ or } \qquad {\rm (b)} \ \  y_+ > y'.
\end{align}
Then we have
\begin{align} \label{eq: La computing formula}
\La( M^\ii[x,y], M^\ii[x',y']) = \sum_{ u \in [x,y]_\ii  , \ v \in [x',y']_\ii  } \la_{u,v}.
\end{align}
\end{proposition}

\begin{proof}
Since the proof is similar, we shall give only the proof of {\rm (a)} case. Let us divide into sub-cases as below.

\noindent
(i) [$x=y > x'_-$] \  If $x>x'$, we have
\begin{align*}
\La(S^\ii_x, M^\ii[x',y']) & = \La(S^\ii_x,M^\ii[x'_+,y'] \hconv S^\ii_{x'}) \\
& \underset{(1)}{=} \La(S^\ii_x,M^\ii[ x'_+,y'])+  \La(S^\ii_x, S^\ii_{x'})= \La(S^\ii_x,M^\ii[x'_+,y'])+ \la^\ii_{x,x'}.
\end{align*}
Here $\underset{(1)}{=}$ holds by \cite[Proposition 2.9]{KKOP22B}  and the fact that $(S^\ii_x,S^\ii_{x'})$ is an unmixed pair.
Then by the induction hypothesis on $|[x',y']_\ii|$, we have
$$
\La(S^\ii_x, M^\ii[x',y']) = \la^\ii_{x,x'} + \sum_{v \in [x'_+,y']_\ii} \la^\ii_{x,v} = \sum_{v \in [x', y' ]_\ii} \la^\ii_{x,v},
$$
as we desired.

Now, the remaining  case of (i) can be described as follows:
$$  x'_-<x =y \le x' \le y'.$$
Since $S^\ii_x$ commutes with $M^\ii[x',y']$ and $M^\ii[x',y'_-]$ by \cite[Proposition 3.27]{KKOP21},
\begin{align*}
\La(S^\ii_x,M^\ii[x',y']) & = - \La(M^\ii[x',y'],S^\ii_x)  = -\La( S^\ii_{y'} \hconv M^\ii[x',y'_-],S^\ii_x)  \\
& = -\La( S^\ii_{y'} ,S^\ii_x) -\La(  M^\ii[x',y'_-],S^\ii_x)  \\
& = (\be_{y'},\be_x) + \La( S^\ii_x, M^\ii[x',y'_-]) = \la^\ii_{x,y'} + \La( S^\ii_x, M^\ii[x',y'_-]),
\end{align*}
Then our assertion follows from the induction hypothesis on $|\,[x',y']_\ii|$.

\noindent
(ii) [$x<y$] \  Assume first that $y>y'$. Then we have
\begin{align*}
\La(M^\ii[x,y],M^\ii[x',y']) & = \La(S^\ii_y \hconv M^\ii[x,y_-],M^\ii[x',y'])  \\
& =  \La(S^\ii_y,M^\ii[x',y']) + \La(M^\ii[x,y_-],M^\ii[x',y']).
\end{align*}
Note that $(S^\ii_y,M^\ii[x',y'])$ is an unmixed pair. Then by induction on $|[x,y]_\ii|$, we have
\eqn
\La(M^\ii[x,y],M^\ii[x',y'])&& = \La(S^\ii_y \hconv M^\ii[x,y_-],M^\ii[x',y'])\\ && = \sum_{v \in [x',y']_\ii} \la^\ii_{y,v} +    \sum_{ u \in [x,y_-]_\ii; \ v \in [x',y']_\ii   } \la^\ii_{u,v}
\eneqn
which yields our assertion for this case.

Now let us assume that $y \le y'$. Then we have
$$x'_- < x < y \le y'.$$
Then for any $u \in [x,y]_\ii$, $S^\ii_u$ commutes with $M^\ii[x',y']$  by \cite[Proposition 3.27]{KKOP21}. By~\cite[Proposition 3.2.13]{KKKO18}, we have
\begin{align*}
\La(M^\ii[x,y],M^\ii[x',y'])  = \sum_{u \in [x,y]_\ii} \La(S^\ii_u,M^\ii[x',y']).
\end{align*}
Then our assertion follows from (i).
\end{proof}

We say that $\ii$-boxes $[a_1,b_1]$ and $[a_2,b_2]$ \emph{commute} 
if we have either
$$
(a_1)_-  < a_2 \le  b_2 < (b_1)_+ \quad \text{ or } \quad (a_2)_-  < a_1 \le  b_1 < (b_2)_+.
$$

The following corollary is proved in \cite[Theorem 4.21]{KKOP2}
in the quantum affine case.
\begin{corollary} \label{cor: commuting Ms}
For commuting $\ii$-boxes  $[a_1,b_1]$ and $[a_2,b_2]$, the modules $M^\ii[a_1,b_1]$ and $M^\ii[a_2,b_2]$ commute.
\end{corollary}

\begin{proof}
By Proposition~\ref{prop: Mi Li}, we have
$$
\La(M^\ii[a_1,b_1],M^\ii[a_2,b_2]) = \sum_{\substack{u \in [a_1,b_1]_\ii \\  v \in [a_2,b_2]_\ii }} \la^\ii_{u,v} = - \La(M^\ii[a_2,b_2],M^\ii[a_1,b_1]),
$$
which implies $\de(M^\ii[a_1,b_1],M^\ii[a_2,b_2])=0$. Thus our assertion follows from Proposition~\ref{prop: simple head and socle}~\eqref{it: de 0 simple}.
\end{proof}

\begin{proposition} \label{prop: Li commute}
For a commuting pair $(M^\ii[x,y],M^\ii[x',y'])$, ~\eqref{eq: La computing formula} holds.
\end{proposition}

\begin{proof}
If the $\ii$-boxes $[x,y]$ and $[x',y']$ satisfies~\eqref{eq: left pre-commutative}, our assertion holds. Thus it is enough to consider when $x \le x'_-$ and $y_+ \le y'$. Since they commute,
$$  \La(M^\ii[x,y],M^\ii[x',y']) =- \La(M^\ii[x',y'],M^\ii[x,y]).$$
If $x'> x_-$ or $y'_+> y$,  Proposition~\ref{prop: Mi Li} tells that
$$  \La(M^\ii[x,y],M^\ii[x',y'])  = - \sum_{  u \in [x,y]_\ii ; \; v \in [x',y']_\ii }  \la^\ii_{v,u} = \sum_{ u \in [x,y]_\ii ; \; v \in [x',y']_\ii } \la^\ii_{u,v},
$$
which implies the assertion. Thus we may assume that $x' \le x_-$. However, in this case, we have
$$
x' \le x_- \le x \le x'_-,
$$
which yields contradiction.
\end{proof}

Let us define the skew symmetric pairing $\rmL_\ii$
on $\Irr(\scrC_w)$ as follows:
\begin{align} \label{eq: def Li}
\rmL_\ii(X,Y) \seteq \sum_{1 \le a,b \le r} (\PBW_\ii(X))_a (\PBW_\ii(Y))_b \; \la^\ii_{a,b}.
\end{align}

The following lemma follows
from Lemma~\ref{lem: PBW addition} and~\eqref{eq: def Li}.

\Lemma \label{cor: Li decomposition}
For $M=\calM^\ii(\bfa)$ with $\bfa \in \Z_{\ge0}^{\K}$, we have
$$    \rmL_\ii(X \hconv M,Y) = \rmL_\ii(X,Y)+ \rmL_\ii(M,Y) \quad \text{ and } \quad \rmL_\ii(X,Y) = -\rmL_\ii(Y,X).$$
\enlemma

\Prop \label{prop: Gr equals L}
For any simple $X,Y$ in $\catC_w$, we have
$$\rmL_\ii(X,Y)=\rmG^R_{\calM^\ii}(X,Y).$$
\enprop
\Proof
Let $\shs$ be the set of simple modules $Y$ in $\catC_w$ such that
$\rmL_\ii(X,Y)=\rmG^R_{\calM^\ii}(X,Y)$ for any simple $X\in\catC_w$, and let
$\shs'$ be the set of simple modules $Y$ in $\catC_w$ such that
$\rmL_\ii(\calM^\ii(\bfa),Y)=\rmG^R_{\calM^\ii}(\calM^\ii(\bfa),Y)$ for any simple
$\bfa\in\Z_{\ge0}^\K$.
By Proposition~\ref{prop: Li commute},
we have
\eqn \rmL_\ii(M^\ii_s,M^\ii_t)=\La(M^\ii_s,M^\ii_t) \quad \text{for any $s,t \in \K$.}\eneqn
 Thus we have
$\calM^\ii(\bfa) \in\shs'$
by Lemma \ref{cor: Li decomposition}.
\medskip
Now, let us show $\shs'\subset\shs$.
Let $Y\in\shs'$. For any simple $X$, there exist
$\bfa,\bfb\in \Z_{\ge0}^\K$ such that $X\hconv \calM^\ii(\bfa)\simeq\calM^\ii(\bfb)$.
Hence we have
\eqn&&
\rmL_\ii(X,Y)+\rmL_\ii(\calM^\ii(\bfa),Y)\underset{(1)}{=}\rmL_\ii(\calM^\ii(\bfb),Y)=
\rmG^R_{\calM^\ii}(\calM^\ii(\bfb),Y)\\
&&\hs{5ex}\underset{(2)}{=}\rmG^R_{\calM^\ii}(X,Y)+\rmG^R_{\calM^\ii}(\calM^\ii(\bfa),Y)=
\rmG^R_{\calM^\ii}(X,Y)+\rmL_\ii(\calM^\ii(\bfa),Y).
\eneqn
Here $\underset{(1)}{=}$ follows from  Lemma~\ref{cor: Li decomposition}
and $\underset{(2)}{=}$
follows from Lemma~\ref{lem:Grconv}.
Hence we have $\rmL_\ii(X,Y)=\rmG^R_{\calM^\ii}(X,Y)$.
Thus we have proved $\shs'\subset\shs$.

Since $\calM^\ii(\bfa)\in\shs'$,
we have $\calM^\ii(\bfa)\in\shs$,
which implies that
$$\rmL_\ii(Y,\calM^\ii(\bfa))=\rmG^R_{\calM^\ii}(Y,\calM^\ii(\bfa))$$ for any simple $Y$.
Hence any simple $Y$ belongs to $\shs'$ and hence to $\shs$.
\QED

\subsection{ Degree and codegree}

In this subsection, we see the relation between $\bfg^R_\calM(X)$ (resp.  $\bfg^L_\calM(X)$) and the degree (resp. codegree) in the (quantum) cluster algebra theory.
For a commuting family $\calM=\{ M_j \mid j \in J \}$ labeled by a finite index set $J$, let us recall the preorder
$\preceq_\calM$ on $\Z_{\ge 0}^{J}$ given in \cite[\S 3.3]{KK19} (see also \cite[Definition 3.1.1]{Qin17}):
\begin{align*}
\bfb' \preceq_\calM \bfb \quad \text{ if and only if } \quad  & {\rm (1)} \; \wt(\calM(\bfb)) =\wt(\calM(\bfb')), \\
& {\rm (2)} \; \La(\calM(\bfb'),M_j) \le  \La(\calM(\bfb),M_j) \text{ for all } j \in J.
\end{align*}

The preorder  $\preceq_\calM$ can be extended to the one on $\Z^{J}$ as follows: For $\bfb,\bfb' \in \Z^{J}$,
$$
\bfb' \preceq_\calM \bfb \text{ if }   \bfb'+\bfa \preceq_\calM \bfb+\bfa \text{ for some } \bfa \in \Z_{\ge 0}^{J}
\text{ such that } \bfb+\bfa,\bfb'+\bfa \in \Z_{\ge 0}^{J}.
$$
We write
$\bfb'\prec_\calM\bfb$,
if $\bfb'\preceq_\calM\bfb$ holds but  $\bfb\preceq_\calM\bfb'$ does not hold. Hence
$\bfb'\prec_\calM\bfb$ if and only if
 $\bfb'\preceq_\calM\bfb$ and there exists $j\in J$ such that
$\La(\calM(\bfb),M_j)<\La(\calM(\bfb'),M_j)$.

\begin{lemma} [cf.{ \cite[Lemma 3.6]{KK19}}] \label{lem: g-vectors and composition series}
Let $X$ be a simple module and $\calM=\{ M_j \mid j \in J\}$ be a \ql family in $\scrC$.
Let  $\bfa \in \Z_{\ge 0}$, $\sfS$, $c(s) \in \Z$ and $\bfb(s) \in \Z_{\ge 0}^J$ $(s \in \sfS)$ be as in \eqref{eq: Laurent}.
Then we have the followings$:$
\bnum
\item \label{it: head} There exists a unique $s_0 \in \sfS$ such that $X \hconv \calM(\bfa) \simeq q^{\,c(s_0)} \calM(\bfb(s_0))$. Moreover, we have
$\bfb(s) \prec_\calM \bfb(s_0)$ for any $s \in \sfS \setminus \{ s_0 \}$.
\item \label{it: socle} There exists a unique $s_1 \in \sfS$ such that $X \sconv \calM(\bfa) \simeq q^{\,c(s_1)} \calM(\bfb(s_1))$. Moreover, we have
$\bfb(s_1) \prec_\calM \bfb(s)$ for any $s \in \sfS \setminus \{ s_1 \}$.
\item \label{it: simplicity} If $s_0 = s_1$, then $\sfS = \{ s_0\}$ and $X \conv \calM(\bfa) \simeq \calM(\bfb(s_0))$.
\item \label{it: lengtht} If  $s_0 \ne s_1$ and there exists no $\bfc \in \Z^\K$ such that $\bfb(s_1) - \bfa \prec_\calM \bfc \prec_\calM \bfb(s_0) - \bfa$, then
$$  [X \conv \calM(\bfa)] =    [X \hconv \calM(\bfa)] +    [X \sconv \calM(\bfa)] \quad \text{ in $K(\scrC_w)$.} $$
\ee
\end{lemma}

\begin{proof}
It follows from Proposition~\ref{prop: simple head and socle} and~\eqref{eq: Laurent}.
\end{proof}

\smallskip

The following proposition is proved for symmetric quiver Hecke algebras and can be extended to general quiver Hecke algebras using the almost same argument:

\begin{proposition}[{\cite[Proposition 3.3]{KK19}}] \label{eq: the order}   For a monoidal cluster $\calM=\{ M_k \mid k \in \K\}$ associated with a quantum seed $(\{ X_k\}_{k\in \K}, L,\tB)$,
$$
\bfb' \preceq_\calM \bfb \quad \text{ if and only if } \quad \bfb-\bfb' = \tB \uv \text{ for some } \uv \in \Z_{\ge0}^\Kex.
$$
In particular, the relation $\preceq_\calM$ is an order on $\Z^\K$.
\end{proposition}

\begin{corollary} \label{cor: pointed}
Let $\calM=\{ M_k \mid k \in \K\}$ be a monoidal cluster associated with a quantum seed $\calS=( \{X_k\}_{k \in \K}, L, B)$.
Then $[X]$ is $\calT(L)$-pointed and $\calT(L)$-copointed for any simple module $X \in \scrC_w$.
\end{corollary}

\begin{remark} For a monoidal cluster $\calM$ associated with a quantum seed $\seed$ and a simple module
$M \in \scrC$,
 the above corollary tells that $\bfg^R_\calM(M)$ and $\bfg^L_\calM(M)$ coincide with the degree and codegree of $[M] \in \calK_\bbA(\scrC) \simeq \scrA_\bbA(\calS)$, respectively.
\end{remark}

\end{document}